\documentclass[a4paper,10pt]{article}
\usepackage{geometry}
\usepackage[utf8]{inputenc}
\usepackage{amssymb}
\usepackage{amsmath}
\usepackage{amsfonts}
\usepackage{amsthm}
\usepackage{lscape}
\usepackage{stmaryrd}
\usepackage{mathrsfs}
\usepackage{wrapfig}
\usepackage{algorithm}
\usepackage{algpseudocode}
\usepackage{amsopn}
\usepackage{mathtools}
\usepackage{pgfplots}
\usepackage{pgfplotstable}
\usepackage{booktabs}
\usepackage{mathabx}
\usepackage{multirow}
\usepackage{cite}
\usepackage[hidelinks,colorlinks=true,allcolors=blue]{hyperref}

\usepackage{tikz}
\usetikzlibrary{matrix,decorations.pathreplacing,calc}

\newcommand{\rank}{\hbox{rank}}
\newcommand{\upd}[1]{{\color{black}#1}} 
\newcommand{\interior}[1]{%
	{\kern0pt#1}^{\mathrm{o}}%
}

\pgfkeys{tikz/mymatrixenv/.style={decoration=brace,every left delimiter/.style={xshift=3pt},every right delimiter/.style={xshift=-3pt}}}
\pgfkeys{tikz/mymatrix/.style={matrix of math nodes,left delimiter=[,right delimiter={]},inner sep=2pt,column sep=1em,row sep=0.5em,nodes={inner sep=0pt}}}
\pgfkeys{tikz/mymatrixbrace/.style={decorate}}

\newcommand{\norm}[1]{\lVert#1\rVert}

\catcode`@=11
\font\tenex=cmex10 
\newdimen\p@renwd
\setbox0=\hbox{\tenex B} \p@renwd=\wd0 
\def\bmat#1{\begingroup \m@th
	\setbox\z@\vbox{\def\cr{\crcr\noalign{\kern2\p@\global\let\cr\endline}}%
		\ialign{$##$\hfil\kern2\p@\kern\p@renwd&\thinspace\hfil$##$\hfil
			&&\quad\hfil$##$\hfil\crcr
			\omit\strut\hfil\crcr\noalign{\kern-\baselineskip}%
			#1\crcr\omit\strut\cr}}%
	\setbox\tw@\vbox{\unvcopy\z@\global\setbox\@ne\lastbox}%
	\setbox\tw@\hbox{\unhbox\@ne\unskip\global\setbox\@ne\lastbox}%
	\setbox\tw@\hbox{$\kern\wd\@ne\kern-\p@renwd\left[\kern-\wd\@ne
		\global\setbox\@ne\vbox{\box\@ne\kern2\p@}%
		\vcenter{\kern-\ht\@ne\unvbox\z@\kern-\baselineskip}\,\right]$}%
	\null\;\vbox{\kern\ht\@ne\box\tw@}\endgroup}
\catcode`@=12

\usepackage{natbib}
\setcitestyle{square}
\setcitestyle{numbers}
\setcitestyle{comma}
\usepackage{caption}
\usepackage{subcaption}

\newcommand{\R}{\mathbb {R}}
\newcommand{\C}{\mathbb {C}}

\newcommand{\diag}{\mathrm{diag}}
\newtheorem{theorem}{Theorem}

\newtheorem{remark}[theorem]{Remark}

\theoremstyle{definition}
\newtheorem{example}[theorem]{Example}

\newcommand{\subop}{\tilde r_j}

\title{Improved ParaDiag via low-rank updates and interpolation}

\author{Daniel Kressner\footnote{MATH-ANCHP, École Polytechnique Fédérale de Lausanne, Station 8, 1015 Lausanne, Switzerland. E-mail: daniel.kressner@epfl.ch} \and Stefano Massei\footnote{Department of Mathematics, University of Pisa, Pisa, Italy. E-mail: stefano.massei@unipi.it} \and Junli Zhu\footnote{School of Mathematics and Statistics, Lanzhou University, Lanzhou 730000, Gansu, PR China. The work of the third author has been supported by National Natural Science Foundation of China (Grant Nos. 11471150, 12161030) and China Scholarship Council  (Contract No. 202006180060). E-mail: zhujl18@lzu.edu.cn}}

\date{}

\pgfplotsset{compat=1.16}
\begin{document}
	
	\maketitle
	
	\begin{abstract}
		This work is concerned with linear matrix equations that arise from the space-time discretization of time-dependent linear partial differential equations (PDEs). Such matrix equations have been considered, for example, in the context of parallel-in-time integration leading to a class of algorithms called ParaDiag. We develop and analyze two novel approaches for the numerical solution of such equations. Our first approach is based on the observation that the modification of these equations performed by ParaDiag in order to solve them in parallel has low rank. Building upon previous work on low-rank updates of matrix equations, this allows us to make use of tensorized Krylov subspace methods to account for the modification. Our second approach is based on interpolating the solution of the matrix equation from the solutions of several modifications. Both approaches avoid the use of iterative refinement needed by ParaDiag and related space-time approaches in order to attain good accuracy. In turn, our new algorithms have the potential to outperform, sometimes significantly, existing methods. This potential is demonstrated for several different types of PDEs. 
	\end{abstract}
	
	
	
	\section{Introduction}
	This paper is concerned with the efficient numerical solution of matrix equations that arise from implicit time integration of large systems of ordinary differential equations (ODEs). More specifically, we treat generalized Sylvester equations of the form
	\begin{equation}\label{eq:main-eq}
		AXB_2^T + MXB_1^T=F,
	\end{equation}
	where $F\in\mathbb R^{n\times n_t}$ and the columns of the unknown matrix $X\in\mathbb R^{n\times n_t}$ contain the approximate solution of the ODE at $n_t$ time steps with constant step size. The matrices $A,M\in\mathbb R^{n\times n}$ are large and sparse, as they arise, for example, from the finite difference/element spatial discretization of a time-dependent partial differential equation (PDE). The matrices $B_1,B_2\in\mathbb R^{n_t\times n_t}$ are banded Toeplitz matrices and contain the coefficients of the time integrator.  We assume throughout this work that $M$ and $B_2$ are invertible. 
	
	Vectorizing~\eqref{eq:main-eq} and using the Kronecker product $\otimes$, we obtain the equivalent $nn_t\times nn_t$ linear system
	\begin{equation}\label{eq:all-at-once-lin}
		\left(B_2\otimes A +B_1\otimes M\right)\mathbf x=\mathbf f,
	\end{equation}
	where $\mathbf x = \mathrm{vec}(X)$ and $\mathbf f = \mathrm{vec}(F)$.
	While the theoretical properties of generalized Sylvester equations are well understood  and various numerical solution strategies have been developed~\cite{benner2013,simoncini2016}, the structure of the matrix coefficients in~\eqref{eq:main-eq} comes with particular challenges that are addressed in this work.
	
	\begin{paragraph}{Illustrative example.} To illustrate~\eqref{eq:main-eq} and motivate our developments, we 
		consider the ODE
		\begin{equation}\label{eq:ode}\begin{cases}M\dot{\mathbf  u}(t)+A\mathbf u(t)= \mathbf f(t),&t\in (0, T],\\
				\mathbf u(0)=\mathbf x_0,\end{cases}
		\end{equation}
		with $\mathbf u(t)\in\mathbb R^n$. When such an ODE arises from the finite element discretization of a linear time-dependent PDE then $A$ and $M$ correspond to the stiffness matrix and the mass matrix, respectively, and $M$ is symmetric positive definite. When using a finite difference discretization, we have $M = I_n$.
		The implicit Euler method with step size $\Delta t = T / n_t$ applied to~\eqref{eq:ode} requires the successive solution of the linear systems
		\begin{equation}\label{eq:time-step}
			(\Delta t^{-1}M+ A)\mathbf x_k =  \mathbf f_k + \Delta t^{-1}M \mathbf x_{k-1},  \qquad k=1,\dots, n_t,
		\end{equation}
		where $\mathbf x_k \approx \mathbf u(k\cdot \Delta t)$.
		Stacking these equations yields an $nn_t\times nn_t$ linear system of the form
		$$
		\begin{bmatrix}
			\Delta t^{-1}M+ A\\
			-\Delta t^{-1}M&\Delta t^{-1}M+ A\\
			&\ddots&\ddots\\
			&&-\Delta t^{-1}M&\Delta t^{-1}M+ A
		\end{bmatrix}\begin{bmatrix}
			\mathbf x_1\\
			\mathbf x_2\\
			\vdots\\
			\mathbf x_{n_t}
		\end{bmatrix}=\begin{bmatrix}\mathbf f_1+ \Delta t^{-1}M\mathbf x_0\\
			\mathbf f_2\\
			\vdots\\
			\mathbf f_{n_t}
		\end{bmatrix}.
		$$
		The system matrix can be rewritten in the form $B_2\otimes A +B_1\otimes M$ from~\eqref{eq:all-at-once-lin} by setting
		\begin{equation} \label{eq:b1b2}
			B_1= \frac 1{\Delta t}\begin{bmatrix}
				1&\\
				-1&1\\
				&\ddots&\ddots\\
				&&-1&1
			\end{bmatrix},\qquad B_2=I_{n_t}.
		\end{equation}
		In turn, the linear system can be rephrased as a matrix equation of the form~\eqref{eq:main-eq}.
		Taking a closer look at \eqref{eq:all-at-once-lin}, we identify the factors that have a role in determining such a peculiar structure. The Kronecker structure of the coefficient matrix is due to the time-independence of the coefficients of the linear ODE~\eqref{eq:ode}. The matrices $B_1,B_2$ are Toeplitz because of the use of a  constant step size. Indeed, this property is maintained when replacing implicit Euler by implicit multistep methods~\cite{Hairer2010} such as BDF (backward differentiation formula).  When using Runge--Kutta integrators~\cite{Hairer2010}, we still get banded Toeplitz matrices $B_1,B_2$ but -- as we will see in  Section~\ref{sec:runge-kutta} -- we have to embed $A, M$ into larger matrices in order to arrive at a matrix equation of the form~\eqref{eq:main-eq}.
	\end{paragraph}
	
	\begin{paragraph}{Existing work.}
		Classical time integration approaches for ODEs proceed by solving the discretized equations, like~\eqref{eq:ode}, sequentially first for $k =1$, then for $k = 2$, and so on. This limits the potential for parallelizing the implementation of the integrator.  The Parareal algorithm~\cite{lions2001,gander2007} avoids this limitation by combining, iteratively, coarse (cheap) time steps that are performed sequentially with fine (expensive) time steps that are performed in parallel. There have been many developments and modifications of Parareal during the last two decades; we refer to~\cite{Gander2015,gander2022} for overviews and a historical perspective. Our work relates to the so called ParaDiag algorithms~\cite{gander2020paradiag}, which proceed somewhat differently from Parareal by explicitly exploiting the structure of the matrix equation~\eqref{eq:main-eq}. \upd{In passing, we mention that linear ODEs with constant coefficients, like the model problem~\eqref{eq:ode}, can be turned into an independent set of linear systems by combining an integral transform, like the Laplace transform, with  quadrature for the numerical evaluation of the inverse transform; see~\cite[Sec. 5.5]{Gander2015} and the references therein.}
		
		The basic idea of ParaDiag is simple:
		If $B_2^{-1} B_1$ can be diagonalized by a similarity transformation then~\eqref{eq:main-eq} decouples into $n_t$ linear systems, which can be solved in parallel. The pitfall of this approach is that the most straightforward choice of time steps, constant time steps, leads to non-diagonalizable matrices; indeed, the matrix $B_2^{-1} B_1$ in~\eqref{eq:b1b2} is one big Jordan block. When choosing a geometrically increasing sequence of time steps, such as $\Delta t_j= \Delta t_1 \tau^{j-1}$ for some $\tau>1$ then $B_2^{-1}B_1$ has mutually distinct eigenvalues and thus becomes diagonalizable~\cite{maday2008}. However, \upd{even when using a second order method in time, e.g., the Crank-Nicolson method,} the condition number of the transformation matrix grows rapidly as $n_t$ increases, which severely limits the number of time steps~\cite{gander2016,gander2019}.
		
		A second class of ParaDiag algorithms allows for uniform time steps and perturbs $B_1, B_2$ to enforce diagonalizability. More specifically, the perturbed system matrix takes the form
		\begin{equation} \label{eq:palpha}
			P_{\alpha}=C_2^{(\alpha)}\otimes A+C_1^{(\alpha)}\otimes M,
		\end{equation}
		where $C_1^{(\alpha)},C_2^{(\alpha)}$ are Strang-type~\cite{Strang1986,bertaccini2002spectrum,noschese2011} $\alpha$-circulant matrices\footnote{On a continuous-time level this strategy can viewed as turning the initial value problem into a boundary value problem with periodic boundary conditions and, hence, it is sometimes called \emph{waveform relaxation}~\cite{wu2018,gander2019b}.} 
		constructed from the Toeplitz matrices $B_1,B_2$. This matrix can be used within a stationary iteration $P_{\alpha}\mathbf x^{(j)}=(P_{\alpha}-B_2\otimes A - B_1\otimes M)\mathbf x^{(j-1)}+\mathbf f$ \cite{liu2020} or as a preconditioner for GMRES~\cite{bertaccini2003block,mcdonald2018} applied to~\eqref{eq:all-at-once-lin}. 
		In the solution of linear systems with $P_{\alpha}$ one takes advantage of the fact that the  matrices $C_1^{(\alpha)},C_2^{(\alpha)}$ are simultaneously diagonalized by a scaled  Fourier transform: 
		$$
		\Omega^*D_{\alpha}C_1^{(\alpha)}D_{\alpha}^{-1}\Omega  =  \diag(\lambda_{1,1},\dots, \lambda_{1,n_t}), \qquad \Omega^*D_{\alpha}C_2^{(\alpha)}D_{\alpha}^{-1}\Omega= \diag(\lambda_{2,1},\dots, \lambda_{2,n_t}),
		$$
		where $\Omega\in\mathbb C^{n_t\times n_t}$ is the discrete Fourier transform matrix and $D_{\alpha}=\diag(1, \alpha^{\frac 1{n_t}},\dots, \alpha^{\frac{n_t-1}{n_t}})$.
		In turn, $
		P_{\alpha}$
		is  block diagonalized:
		$$(\Omega^*D_{\alpha}\otimes I_n)P_{\alpha}(D_{\alpha}^{-1}\Omega \otimes I_n)= \begin{bmatrix}\lambda_{1,1}M+\lambda_{2,1}A\\
			&\ddots\\
			&&\lambda_{1,n_t}M+\lambda_{2,n_t}A\end{bmatrix}.$$
		Using the fast Fourier transform (FFT), this allows us to compute the solution of a linear system, with $P_{\alpha}$ by combining FFT with the (embarrassingly parallel) solution of $n_t$ linear systems with the matrices $(\lambda_{1,j}M +\lambda_{2,j}A$, $j = 1,\ldots,n_t$.
		This procedure, rephrased for the matrix equation corresponding to $P_{\alpha}\mathbf x=\mathbf f$, is summarized in Algorithm~\ref{alg:fast-diag}. \upd{Note that, even though the original ODE~\eqref{eq:ode} features real matrices, Algorithm~\ref{alg:fast-diag} requires the solution of \emph{complex} linear systems, which is more expensive (often by a factor $2$--$4$) than solving real linear systems. This disadvantage of having to pass to complex arithmetic is shared by all diagonalization based methods discussed in this work.} 
		\begin{algorithm}[H] 
			\small 
			\caption{Solution of $AX(C_2^{(\alpha)})^T+MX(C_1^{(\alpha)})^T=F$}\label{alg:fast-diag}
			\begin{algorithmic}[1]
				\Procedure{fast\_diag\_solve}{$A,M,B_1,B_2, F,\alpha$}
				\State Compute eigenvalues $\lambda_{1,1},\ldots,\lambda_{1,n_t},\lambda_{2,1},\ldots,\lambda_{2,n_t}$ by FFT.
				\State Compute $X \gets F D_{\alpha}$ by column scaling and $X\gets X\Omega^*$ by FFT.   
				\For{$j=1,\dots, n_t$}
				\State Solve linear system $X(:, j)\gets (\lambda_{1,j}M + \lambda_{2,j}A)^{-1}X(:, j)$
				\EndFor 
				\State Compute $X \gets X \Omega$ by FFT and $X \gets X D_{\alpha}^{-1}$ by column scaling. 
				\EndProcedure
			\end{algorithmic}
		\end{algorithm} 
		The use of the preconditioner $P_\alpha$ can dramatically accelerate the convergence of an iterative solver for~\eqref{eq:all-at-once-lin}. This has been observed for parabolic problems in~\cite{mcdonald2018} when using  $P_{\alpha}$ with $\alpha=\pm 1$ as a preconditioner in GMRES. In~\cite{gander2020paradiag}, wave propagation problems (integrated in time with a leap-frog finite difference scheme) have been successfully treated by choosing a relatively small value for $0<\alpha<1$. Note, however,  that $\alpha$ cannot be chosen arbitrarily small because otherwise the condition number of the matrix $D_\alpha$ (given by $\alpha^{(1-n_t)/n_t} \approx 1/\alpha$) explodes and, in turn, roundoff error impedes convergence. In practice, this means that the choice of $\alpha$ needs to strike a compromise and there is no choice of $\alpha$ that yields sudden convergence. In turn, several iterations are needed to attain good accuracy, which adds significant computational overhead and communication cost because every iteration calls Algorithm~\ref{alg:fast-diag} once. 
	\end{paragraph}
	
	\begin{paragraph}{New contributions.}
		In this paper, we propose two novel strategies that avoid the overhead incurred by iterative procedures for solving~\eqref{eq:all-at-once-lin}. 
		
		Our first method, presented in Section~\ref{sec:low-rank}, exploits the trivial observation that $C_1^{(\alpha)} - B_1$ and $C_2^{(\alpha)} - B_2$ have very low rank. Thus, the Sylvester equation addressed by Algorithm~\ref{alg:fast-diag} can be viewed as a low-rank modification of the original equation~\eqref{eq:main-eq}, which allows us to apply the low-rank update  from~\cite{kressner2018} for linear matrix equations. This update corrects the output of Algorithm~\ref{alg:fast-diag} by solving a Sylvester equation with the same coefficients as~\eqref{eq:main-eq} but with a different right-hand side that has low rank. The main novelty of Section~\ref{sec:low-rank} is in the analysis of the low-rank structure of the solution to this equation and the development of efficient algorithms.
		
		Our second method, presented in Section~\ref{sec:ev-int}, takes an interpolation point of view: Whenever we solve the linear system $P_{\alpha}\mathbf x=\mathbf f$ we are actually evaluating a vector-valued function $\mathbf x(\alpha)$, depending on the complex parameter $\alpha$, such that $\mathbf x(0)$ corresponds to the solution of \eqref{eq:all-at-once-lin}. More precisely, this function can be efficiently and accurately evaluated on a scaled unit circle. We suggest to construct an approximation for $\mathbf x(0)$ via a linear combination of $\mathbf x(\alpha)$ evaluated at all $d$th order roots of a complex number for some small value of $d$. This strategy is highly parallelizable as it requires to solve $d$ completely independent linear systems with coefficient matrices of the form $P_{\alpha}$. Moreover, it applies also to the cases where the matrices $\delta B_j^{(\alpha)}$ are not low-rank, for example, when a time fractional differential operator is involved.  
		
		Section~\ref{sec:runge-kutta} extends our framework to implicit Runge--Kutta time integrators.
		
		Section~\ref{sec:experiments} is dedicated to a range of numerical experiments that highlight the advantages of our method for a broad range of situations.
	\end{paragraph}
	
	\section{Low-rank update approach}\label{sec:low-rank}
	
	For a lower triangular banded Toeplitz matrix $B$, Strang's $\alpha$-circulant preconditioner \cite{noschese2011} is obtained from reflecting the lower band to the top right corner and  multiplying it by $\alpha$: \upd{
		$$
		\left[\begin{smallmatrix}
			b_0\\
			b_1&b_0\\
			\vdots&\ddots&\ddots\\
			b_{w}&\ddots&\ddots&\ddots\\
			&\ddots&\ddots&\ddots&\ddots\\
			&&b_w&\ldots&b_1&b_0
		\end{smallmatrix}\right]\quad\longrightarrow\quad \left[\begin{smallmatrix}
			b_0\\
			b_1&b_0\\
			\vdots&\ddots&\ddots\\
			b_{w}&\ddots&\ddots&\ddots\\
			&\ddots&\ddots&\ddots&\ddots\\
			&&b_w&\ldots&b_1&b_0
		\end{smallmatrix}\right] +\alpha\left[\begin{smallmatrix}
			\phantom{\dots}&\phantom{\dots}&b_w&\ldots&b_1\\
			&&&\ddots&\vdots\\
			&&&&b_{w}\\
			&\phantom{vdots}\\
			&\phantom{vdots}\\
			&\phantom{vdots}\\
			&\phantom{vdots}\\
			&\phantom{vdots}\\
			&\phantom{vdots}\\
			&\phantom{vdots}\\
		\end{smallmatrix}\right].
		$$}
	In particular, this implies that the rank of this modification is bounded by the width of the lower bandwidth. Applied to the time stepping matrices $B_1$, $B_2$, we obtain that
	\[
	C_1^{(\alpha)}=B_1+\delta B_1^{(\alpha)}, \qquad C_2^{(\alpha)}=B_2+\delta B_2^{(\alpha)}
	\]
	for certain low-rank matrices $\delta B_1^{(\alpha)},\delta B_2^{(\alpha)}$. For instance, for the implicit Euler method (see~\eqref{eq:b1b2}) we have $\delta B_1^{(\alpha)}=-\frac{\alpha}{\Delta t}\mathbf e_1 \mathbf e_{n_t}^T$, where $\mathbf e_j$ denoting the $j$th unit vector of length $n_t$, and $\delta B_2^{(\alpha)}=0$. 
	Inspired by~\cite{kressner2018} we write the solution of \eqref{eq:main-eq} as $X=X_0+\delta X$, where 
	\begin{align}
		AX_0(C_2^{(\alpha)})^T+MX_0(C_1^{(\alpha)})^T &= F,\label{eq:circulant-eq}\\
		A\delta XB_2^T+M\delta XB_1^T &=AX_0(\delta B_2^{(\alpha)})^T + MX_0(\delta B_1^{(\alpha)})^T.\label{eq:update-eq}
	\end{align}
	After equation~\eqref{eq:circulant-eq} is solved via Algorithm~\ref{alg:fast-diag}, the correction equation~\eqref{eq:update-eq} becomes a generalized Sylvester equation in $\delta X$ that has the same coefficients as~\eqref{eq:main-eq} but a different right-hand side of rank at most $\rank(\delta B_1^{(\alpha)})+\rank(\delta B_2^{(\alpha)})$. This remarkable property enables us to make use of low-rank solvers for linear matrix equations (see~\cite{simoncini2016} for an overview) in order to approximate $\delta X$. Note that, the choice of $\alpha\neq 0$ has no influence on this strategy; for this reason we always choose $\alpha=1$ when using the low-rank update approach. The resulting procedure is summarized in Algorithm~\ref{alg:update}. 
	\begin{algorithm}[H] 
		\small 
		\caption{}\label{alg:update}
		\begin{algorithmic}[1]
			\Procedure{low\_rank\_update}{$A,M,B_1,B_2, F$}
			\State $X_0\gets \textsc{fast\_diag\_solve}(A,M,B_1,B_2, F,1)$\label{step:X0}
			\State Compute a low-rank factorization $UV^*=AX_0(\delta B_2^{(1)})^T + MX_0(\delta B_1^{(1)})^T$\label{step:UV}
			\State $\delta X\gets \textsc{low\_rank\_solver}(A,M,B_1,B_2,U,V)$\label{step:rksm} \Comment{Call to a low-rank solver for matrix equations}
			\State\Return $X_0+\delta X$
			\EndProcedure
		\end{algorithmic}
	\end{algorithm} 
	
	\begin{remark}
		It is not uncommon that the right-hand side matrix $F$ of~\eqref{eq:main-eq} has (numerically) low rank because,  for example, the inhomogeneity $\mathbf f$ of the ODE~\eqref{eq:ode} is constant or varies smoothly with respect to $t$. In such a situation, a low-rank solver can be applied directly to~\eqref{eq:main-eq} and the following approach has been suggested by Palitta~\cite{palitta2021}: A block Krylov subspace method is used for reducing the spatial dimension combined with an application of the Shermann--Morrison--Woodbury formula to the resulting reduced equation. Our low-rank approach has the advantage that it is agnostic to properties of $F$, at the expense of having to solve~\eqref{eq:circulant-eq}. In principle, Palitta's method can be applied to solve the correction equation~\eqref{eq:update-eq} but we found it more efficient for large $n_t$  to use a two-sided approach that reduces the temporal dimension as well.
	\end{remark}
	
	\subsection{Solution of correction equation} \label{sec:solcorrection}
	
	In the following, we discuss the application of low-rank solvers to the correction equation~\eqref{eq:update-eq}. In view of the (assumed) invertibility of $M$ and $B_2$ we replace \eqref{eq:update-eq} by the following Sylvester matrix equation:
	\begin{equation}\label{eq:lr-rhs}
		\widetilde A\delta X+\delta X\widetilde B^T=C,
	\end{equation}
	where \[\widetilde A= M^{-1}A, \quad \widetilde B= B_2^{-1}B_1, \quad C= M^{-1}(AX_0(\delta B_2^{(1)})^T + MX_0(\delta B_1^{(1)})^T)B_2^{-T}\] and $\rank(C)\le k:=\rank(\delta B_1^{(1)})+\rank(\delta B_2^{(1)})\ll \min\{n,n_t\}$.
	For given factorizations $\delta B_1^{(1)}=U_1V_1^*$ and $\delta B_2^{(1)}=U_2V_2^*$, a low-rank factorization $C=UV^*$, with $U\in\mathbb R^{n\times k}$ and $V\in\mathbb R^{n_t\times k}$, is obtained by setting
	\begin{equation} \label{eq:facC}
		U = [M^{-1}AX_0 V_2, \  X_0V_1], \qquad V=[B_2^{-1}U_2,\ B_2^{-1}U_1].
	\end{equation}
	
	Alternatively, when $M$ is symmetric positive definite and its Cholesky factorization $M=LL^T$ is available, a possible symmetry of $A$ is preserved by considering
	\begin{equation}\label{eq:update-chol}
		\widecheck A\delta \widecheck{X}+ \delta \widecheck{X}\widecheck B^T= \widecheck C,
	\end{equation}
	with $\widecheck{A}= L^{-1}AL^{-T}$, $\widecheck B=\widetilde B$, $\widecheck C= L^{-1}(AX_0(\delta B_2^{(1)})^T + MX_0(\delta B_1^{(1)})^T)B_2^{-T}$, and $\delta\widecheck X= L^{T}\delta X$. In analogy to~\eqref{eq:facC} an explicit rank-$k$ factorization of $\widecheck C$ can be derived. In the rest of this section we focus on analyzing and solving \eqref{eq:lr-rhs} but we remark that our findings extend to  \eqref{eq:update-chol} with minor modifications.
	
	Various numerical methods have been developed to treat large-scale Sylvester equations with low-rank right-hand side, including the \emph{Alternating Direction Implicit  method} (ADI)~\cite{ellner1986,peaceman1955} and the \emph{Rational Krylov Subspace Method} (RKSM)~\cite{benner2009,grimme1997krylov,ruhe1984}. These methods compute approximate solutions of the form $\widetilde{\delta X}= WYZ^*$, where $W$ and $Z$ are bases of (block) rational Krylov subspaces generated with the matrices $\widetilde A,\widetilde B$ and starting block vectors $U$ and $V$, respectively. More precisely, given two families $\mathbf \xi=\{\xi_1,\dots,\xi_{\ell}\}$, $\mathbf \psi=\{\psi_1,\dots,\psi_{\ell}\}\subset \mathbb C$ of so-called \emph{shift parameters}, $W$ and $Z$ are chosen as bases of the following subspaces:
	\begin{equation} \label{eq:defrk}
		\begin{array}{ll}
			\mathcal{RK}(\widetilde A,U,\mathbf \xi)&= \mathrm{span}\big((\xi_1I-\widetilde A)^{-1}U,\dots, (\xi_{\ell}I-\widetilde A)^{-1}U\big),\\[0.05cm]
			\mathcal{RK}(\widetilde B,V,\mathbf \psi)&= \mathrm{span}\big((\psi_1I-\widetilde B)^{-1}V,\dots, (\psi_{\ell}I-\widetilde B)^{-1}V\big).
		\end{array}
	\end{equation}
	When shifts are repeated, the matrix power is increased. For example, when a single shift $\xi_1$ is repeated $\ell$ times, one uses
	$\mathrm{span}\big((\xi_1I-\widetilde A)^{-1}U,\dots, (\xi_{1}I-\widetilde A)^{-\ell}U\big)$.
	We refer to the relevant literature \cite{benner2009,druskin2011,simoncini2016} for a complete description of the computation of the factors $W,Y,Z$ in these methods.
	
	Let us comment on some implementation aspects of RKSM, which will be the method of choice in our implementation of Algorithm~\ref{alg:update}, line~\ref{step:rksm}. 
	In RKSM, the middle factor $Y$ in the approximation $\widetilde{\delta X}= WYZ^*$ is obtained by solving a compressed matrix equation, which does not contribute significantly to the overall computational time as long as the dimensions of the Krylov subspaces remain modest. Computing $W$ and $Z$ requires to solve shifted linear systems with the matrices $\widetilde A,\widetilde B$; these operations can be performed without forming $M^{-1}A$ and $B_2^{-1}B_1$ explicitly, by means of the relations
	\begin{align*}
		(\xi_jI-\widetilde A)^{-1}\mathbf v= (\xi_jM-A)^{-1}M\mathbf v,\qquad 
		(\psi_jI-\widetilde B)^{-1}\mathbf v= (\psi_jB_2-B_1)^{-1}B_2\mathbf v.
	\end{align*}
	In particular, we can still leverage properties like sparsity in $A$ and $M$, or the bandedness of $B_1,B_2$, when  solving shifted linear systems with $\widetilde A,\widetilde B$. Our implementation of RKSM relies on the Matlab toolbox \texttt{rk-toolbox}~\cite{berljafa2015} for generating the factors $W,Z$.
	Finally and importantly, the selection of the shift parameters $\mathbf \xi,\mathbf \psi$ strongly affects the convergence of RKSM and we discuss suitable choices in the following section.
	
	\subsection{Low-rank approximability of the correction equation} \label{sec:zolotarev}
	
	In this section, we analyze the numerical low-rank structure of the solution $\delta X$ to the correction equation~\eqref{eq:lr-rhs}, which also yields suitable choices of  shift parameters in RKSM \upd{discussed above}. 
	
	\upd{By the Eckart–Young–Mirsky theorem, $\delta X$ admits good low-rank approximations
		if its singular values $\sigma_j(\delta X)$ decay rapidly. Several works, including~\cite{antoulas2002,baker2015,penzl2000,sabino2007}, have established singular value decay bounds for solutions of Sylvester matrix equations.
		In particular, the framework presented in \cite{beckermann2019} applies to~\eqref{eq:lr-rhs} and shows that~\cite[P. 323]{beckermann2019}
		\begin{equation} \label{eq:beckermanntownsend}
			\sigma_{1+kj}(\delta X) \le \min_{r(z)\in\mathcal R_{j,j}} \|r(\widetilde A)\|_2 \|r(-\widetilde B)^{-1}\|_2  \|\delta X\|_2,
		\end{equation}
		where $\mathcal R_{j,j}$ denotes the set of rational functions having numerator and denominator of degree at most $j$. Choosing a good candidate for the rational function $r$ is difficult, especially when $\widetilde A$ and/or $\widetilde B$ are non-normal. To circumvent this difficulty, we loosen the bound~\eqref{eq:beckermanntownsend} by considering the numerical range
		$\mathcal W(\widetilde A) = \{ x^* A x\colon \|x\|_2 = 1\} \subset \mathbb C$. A  result by Crouzeix and Palencia~\cite{crouzeix2017} states that
		\[
		\|r(\widetilde A)\|_2 \le (1+\sqrt 2) \max_{z \in \mathcal W(\widetilde A)} 
		|r(z)| \le (1+\sqrt 2) \max_{z \in E} 
		|r(z)|,
		\]
		where the constant $1+\sqrt 2$ can be omitted when $A$ is normal and the set $E\supseteq\mathcal W(\widetilde A)$ is chosen to feature a simple geometry, for which the eventual rational approximation problem admits explicit solutions.
		Similarly, one has 
		\[
		\|r(-\widetilde B)^{-1}\|_2 \le 
		(1+\sqrt 2) \max_{z \in F} |r(z)|^{-1} = 
		(1+\sqrt 2) \big( \min_{z \in F} |r(z)| \big)^{-1}
		\]
		for any set $F \subset \mathbb C$ with $F\supseteq \mathcal W(-\widetilde B)$.
		Inserting these inequalities into~\eqref{eq:beckermanntownsend}, we arrive at
		\[
		\frac{\sigma_{1+kj}(\delta X)}{\norm{\delta X}_2}\leq (1+\sqrt 2)^a Z_j(E, F), \qquad Z_j(E, F):=\min_{r(z)\in\mathcal R_{j,j}}\frac{\max_{z\in E} |r(z)|}{\min_{z\in F} |r(z)|},
		\]
		where the exponent $a$ satisfies 
	}%
	$a = 0$ if both $\widetilde A,\widetilde B$ are normal, $a = 1$ if one of the two matrices is normal, and $a = 2$ otherwise.

	The quantity $Z_j(E,F)$ and the related optimization problem are known in the literature as the \emph{Zolotarev number} and \emph{Zolotarev problem} for the sets $E$ and $F$ \cite{zolotarev1877,beckermann2019}. Studying the Zolotarev problem is not only important from a theoretical but also from a computational point of view because
	identifying an extremal rational function $r_j^*(z)$ for $Z_j(E,F)$ allows for the fast solution of \eqref{eq:lr-rhs}. Indeed, running $j$ steps of ADI or RKSM with shift parameters given
	by the zeros and poles of $r_j^*(z)$ ensures a convergence rate not slower than the decay rate of $Z_j(E,F)$~\cite{lebedev1977zolotarev,wachspress2013,beckermann2011}.
	\upd{A major complication of this construction is the need for choosing the sets $E,F$. On the one hand, it is desirable to choose $E,F$ such that they enclose the numerical ranges of $\widetilde A,\widetilde B$ as tightly as possible. On the other hand, 
	}%
	explicit solutions for $Z_j(E,F)$ are known only for a few selected geometries of $E$ and $F$, including: two disjoint real intervals~\cite{zolotarev1877}, the two connected components of the complement of an annulus~\cite{gonchar1969}, two disjoint discs~\cite{starke1992}. 
	Therefore, one usually encloses the numerical range by a disc or an interval (in the case of Hermitian matrices).
	In Sections~\ref{sec:twodiscs} and~\ref{sec:discint} below, we discuss two geometries in more detail: the case of two discs and the case of a disc and an interval. These enclosures will be used for generating the shift parameters in Algorithm~\ref{alg:update} for most of our numerical examples.
	
	\subsubsection{Solution of the Zolotarev problem for two disjoint discs} \label{sec:twodiscs}
	
	Letting $\mathbb B(x,\rho)$ denote the closed disc with center $x\in\mathbb C$ and radius $\rho>0$, we consider the Zolotarev problem for
	$$E=\mathbb B(x_E, \rho_E),\qquad  F=\mathbb B(x_F, \rho_F),\qquad |x_E-x_F|> \rho_E+\rho_F.$$  
	Starke~\cite{starke1992} addressed the case $x_E, x_F \in \R$ by explicitly constructing a rational function $r_1^*(z)=\frac{z-q^*}{z-p^*}$ with $q^*\in \interior{E}$, $p^*\in \interior{F}$, and $|r_1^{*}(z)|$ being constant on $\partial E$ and on $\partial F$. In view of the (generalized) \emph{near-circularity criterion} \cite[Theorem 3.1]{starke1992}, this implies the optimality of the rational function $[r_1^*(z)]^j$ for $Z_j(E,F)$, $j\geq 1$.
	
\begin{figure}
		\centering
		\includegraphics[width=\textwidth]{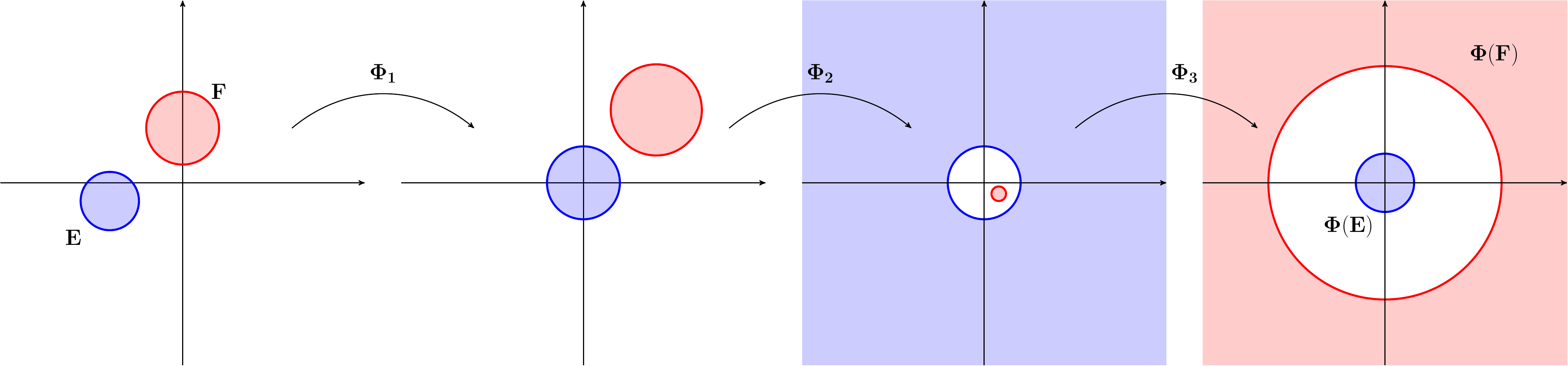}
		\caption{Action of the Moebius map $\Phi=\Phi_3\circ\Phi_2\circ\Phi_1$ on the sets $E$ (blue) and $F$ (red).}
		\label{fig:moeb}
	\end{figure}

	The general case of two disjoint discs (with $x_E,x_F \in\mathbb C$ not necessarily on the real line)
	can be derived from Starke's result by relying on the invariance property of Zolotarev problems under Moebius transformations; see also~\cite[Example VIII.4.2]{Saff1997}. For completeness and convenience of the reader, we provide the complete transformation that maps $E,F$ to the complement of an annulus. More precisely, a Moebius map $\Phi$ is constructed such that $\Phi(F)$ is a closed disc centered at the origin and $\Phi(E)$ is the complement of a larger open disc centered at the origin.  We compose $\Phi$ from three simpler Moebius maps $\Phi_j$, $j=1,2,3$, as follows (see also  Figure~\ref{fig:moeb}):
	\begin{itemize}
		\item $\Phi_1(z)= \frac{z-x_E}{\rho_E}$ maps $E$ to the unit disc while $F$ remains a disjoint disc,
		\item $\Phi_2(z)= z^{-1}$ maps $\Phi_1(E)$ to the exterior of the unit disc and $\Phi_1(F)$ to a disc enclosed by the unit circle,
		\item $\Phi_3(z)=\frac{z-\alpha}{z-\beta}$ maps $\Phi_2(\Phi_1(F))$ to the exterior of a disc with center $0$ and $\Phi_2(\Phi_1(E))$ to a disc centered at $0$ and enclosed by $\partial\Phi_3(\Phi_2(\Phi_1(F)))$.
	\end{itemize} 
	The parameters $\alpha, \beta \in \C$ are chosen as the so-called \emph{common inverse
		points} for the circles $\partial \Phi_2(\Phi_1(E))=\partial \mathbb B(0,1)$ and $\partial \Phi_2(\Phi_1(F))$ \cite[Section 4.2]{henrici1988}. Setting $ \Phi_2(\Phi_1(F))=:\mathbb B(\widetilde x_F, \widetilde \rho_F)$ the values of $\alpha,\beta$ from satisfying the two equations
	\[
	\alpha\overline\beta=1, \qquad (\alpha - \widetilde x_F)\overline{(\beta - \widetilde x_F)}=\widetilde \rho_F^2.
	\]
	Clearly, one solution is given by 
	\begin{equation}  \label{eq:alphabeta}
		\beta=\alpha^{-1}, \qquad \alpha=\bigg[|\widetilde x_F|^2+1-\widetilde \rho_F^2+\sqrt{(|\widetilde x_F|^2+1-\widetilde \rho_F^2)^2-4|\widetilde x_F|^2}\bigg]/(2\overline{\widetilde x_F}).
	\end{equation}
	Observe that, $\Phi_1(F)$ is a disc with center $\Phi_1(x_F)\neq 0$ and radius $\rho_F/\rho_E$. In particular,  $\Phi_1(F)$ attains its maximum and minimum modulus at the elements
	\[
	x_1=\Phi_1(x_F) + \frac{\Phi_1(x_F)}{|\Phi_1(x_F)|}\frac{\rho_F}{\rho_E},\qquad x_2=\Phi_1(x_F) - \frac{\Phi_1(x_F)}{|\Phi_1(x_F)|}\frac{\rho_F}{\rho_E}.
	\]
	Now, $\Phi_2(z)=z^{-1}$ maps $x_1$ and $x_2$ into the minimum and maximum modulus elements of $\mathbb B(\widetilde x_F, \widetilde \rho_F)$, respectively. This implies that 
	$\widetilde x_F,\widetilde \rho_F$ take the following values:
	\[
	\widetilde x_F=\frac{\Phi_2(x_1) + \Phi_2(x_2)}2,\qquad \widetilde \rho_F=\frac{|\Phi_2(x_1) - \Phi_2(x_2)|}2.
	\]
	Together with~\eqref{eq:alphabeta}, this allows us to compute $\alpha,\beta$ explicitly and evaluate $\Phi$.
	The complete expressions for $\Phi$ and $\Phi^{-1}$ read as follows:
	\[
	\Phi(z)= \frac{\rho_E-\alpha(z-x_E)}{\rho_E-\overline{\alpha}^{-1}(z-x_E)},\qquad \Phi^{-1}(z)= \frac{(\overline{\alpha}^{-1}x_E+\rho_E)z-\alpha x_E-\rho_E}{\overline{\alpha}^{-1}z-\alpha}.
	\]
	
	By the near-circularity criterion\upd{~\cite[Theorem 3.1]{starke1992}}, an extremal rational function for the complement of an annulus centered at $0$ is $z^{j}$, having the only pole $\infty$ and the only zero at $0$. Thus, an optimal rational function for the original problem is given by $r_j^*(z)=\left(\frac{z-\Phi^{-1}(0)}{z-\Phi^{-1}(\infty)}\right)^{j}$.
	The pole $p^*$ and zero $q^*$ are given by
	\begin{equation} \label{eq:singleshift}
		p^*=\Phi^{-1}(\infty)= x_E+\frac{\rho_E}{\overline{\alpha}^{-1}},\qquad q^*=\Phi^{-1}(0)= x_E+\frac{\rho_E}{\alpha}.    
	\end{equation}
	Since $\Phi$ maps the boundaries of $E$ and $F$ to the boundaries of the annulus, it follows that the Zolotarev numbers are given by 
	\begin{equation}\label{eq:decay-zol-2DISCS}
		Z_j(E,F)=\eta^j, \qquad 
		\eta=
		\frac{|\Phi(x_E+\rho_E)|}{|\Phi(x_F+\rho_F)|}=\left| \frac{(1-\alpha)(\rho_E-\overline{\alpha}^{-1}(x_F+\rho_F-x_E))}{(1-\overline{\alpha}^{-1})(\rho_E-\alpha(x_F+\rho_F-x_E))}\right|<1.
	\end{equation}
	
	\subsubsection{Quasi-optimal solution of the Zolotarev problem for a disc and an interval} \label{sec:discint}
	
\begin{figure}
		\centering
		\includegraphics[width=.75\textwidth]{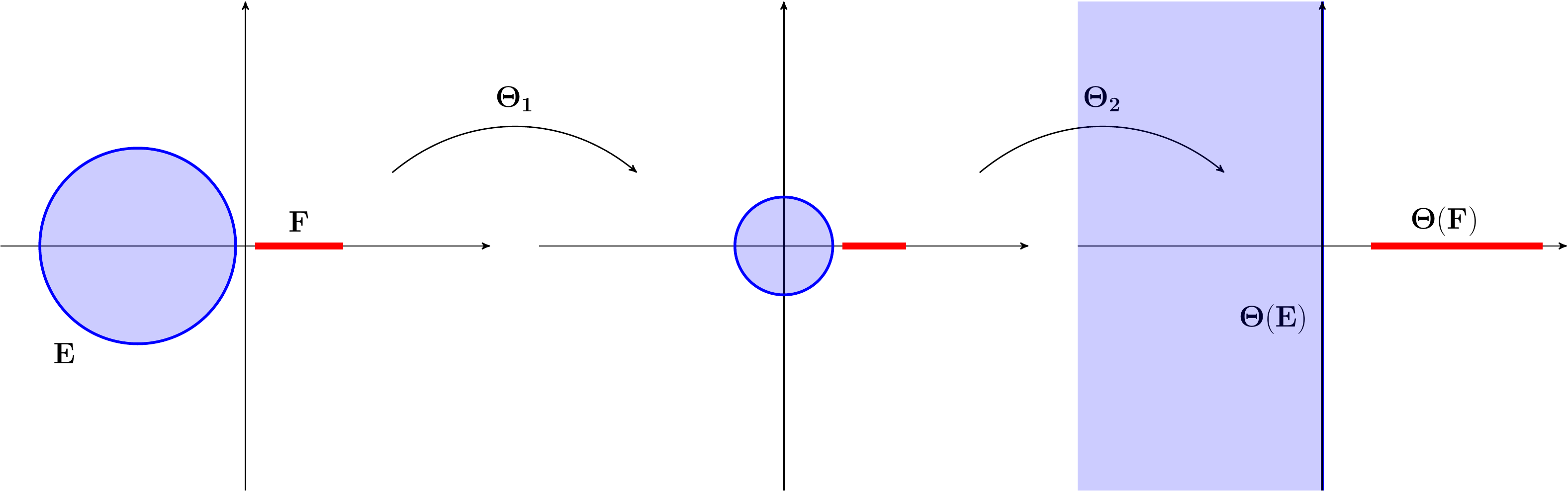}
		\caption{Action of the Moebius map $\Theta=\Theta_2\circ\Theta_1$ on the sets $E$ (blue) and $F$ (red).}
		\label{fig:cayl}
	\end{figure}

	Let us now consider the situation of a disc and a disjoint interval: 
	$$
	E=\mathbb B(x_E,\rho_E),\qquad F =[a,b],\qquad x_E\in \mathbb R,\quad x_E+\rho_E<a\text{ or } x_E-\rho_E>b.
	$$
	Similarly to the previous section, we build a Moebius transformation that recasts $Z_j(E,F)$ as a Zolotarev problem for which a (quasi-optimal) solution is known. More specifically, we let $\Theta=\Theta_2\circ \Theta_1$ where (see also Figure~\ref{fig:cayl}): 
	\begin{itemize}
		\item $\Theta_1(z)=\Phi_1(z)$ maps $E$ to the unit disc while $F$ remains a real disjoint interval,
		\item $\Theta_2(z)=\frac{z+1}{z-1}$ maps $\Theta_1(E)$ to the closed left half of the complex plane and $\Theta_1(F)$ to an interval $[\tilde a, \tilde b]$ of positive real numbers.
	\end{itemize}  
	A straightforward calculation yields explicit expressions for $\Theta$ and its inverse:
	\[
	\Theta(z)=\frac{z-x_E+\rho_E}{z-x_E-\rho_E},\qquad \Theta^{-1}(z)=\frac{(x_E+\rho_E)z-x_E+\rho_E}{z-1}.
	\]%
	By the maximum modulus principle, it suffices to consider the boundary of $\Theta(E)$ and thus
	the transformed Zolotarev problem takes the form
	\begin{equation}\label{eq:zol-imag}
		Z_j(\partial \Theta(E), \Theta(F))=Z_j(\mathbf i\mathbb R, [\tilde a, \tilde b])=\min_{r(z)\in\mathcal R_{j,j}}\frac{\max_{z\in \mathbf i\mathbb R}|r(z)|}{\min_{z\in[\tilde a, \tilde b]}|r(z)|}.
	\end{equation}
	Let $\subop$ denote the extremal rational function for the Zolotarev problem $Z_j( [-\tilde b, -\tilde a], [\tilde a, \tilde b])$, that is, the imaginary axis is replaced by $[-\tilde b, -\tilde a]$. This function has been extensively studied in the literature; see~\cite[Sec. 3.1]{beckermann2019} and the references therein. In particular, it holds that
	$$
	\subop(z)=\prod_{i=1}^j \frac{z+\tilde \psi_{j,i}}{z-\tilde \psi_{j,i}},
	$$
	where the poles $\tilde \psi_{j,i} \in [\tilde a,\tilde b]$, $i=1,\dots, j$, can be expressed in closed form via elliptic integrals. Note that $\subop(-z) = 1/\subop(z)$ and the maximum absolute value of $\subop$ is $1$, which is assumed throughout the imaginary axis. By~\cite[Corollary 3.2]{beckermann2019}, it holds that $Z_j( [-\tilde b, -\tilde a], [\tilde a, \tilde b]) \le 4 \eta^{-2j}$ with $\eta=\exp\left(\dfrac{\pi^2}{2\log(4\tilde b/\tilde a)}\right)$. Inserting $\subop$ into~\eqref{eq:zol-imag} gives the following upper bound for $Z_j(\mathbf i\mathbb R, [\tilde a, \tilde b])$:
	\begin{equation} \label{eq:defdeltaj}
		\delta_j :=  
		\frac{\max_{z\in \mathbf i\mathbb R}|\subop(z)|}{\min_{z\in[\tilde a, \tilde b]}|\subop(z)|} = \frac{1}{\min_{z\in[\tilde a, \tilde b]}|\subop(z)|} = \sqrt{Z_j( [-\tilde b, -\tilde a], [\tilde a, \tilde b])} \le 2 \eta^{-j}.
	\end{equation}
	It has been shown in~\cite{bailly2000,druskin2009} that $\subop$ satisfies the necessary optimality conditions for \eqref{eq:zol-imag} and that it is within a factor $2$ of the optimal solution:  $\delta_j \le 2 Z_j(\mathbf i\mathbb R, [\tilde a, \tilde b])$.
	To the best of our knowledge, it is still an open question whether $\subop$ is, in fact, an extremal rational function for $\eqref{eq:zol-imag}$. 
	In summary, we have
	\begin{equation}\label{eq:decay-zol-DI}
		Z_j(E,F)\leq 2 \exp\left(\frac{-j\pi^2}{2\log(4\tilde b/\tilde a)}\right)
	\end{equation}
	with a quasi-optimal solution for $Z_j(E,F)$ given by \begin{equation}\label{eq:inv-cayl}\widehat r_j(z)=\prod_{i=1}^j\frac{z-\Theta^{-1}(-\tilde \psi_{j,i})}{z-\Theta^{-1}(\tilde \psi_{j,i})}.\end{equation}
	
	The poles of the rational function~\eqref{eq:decay-zol-DI} are not nested, that is, the poles of $\widehat r_j(z)$ are different from the ones of $\widehat r_{j+1}(z)$. This is a disadvantage when using the poles of $\widehat r_j(z)$ in an iterative method (such as ADI or RKSM) where usually the number of steps (the parameter $j$) is not known a priori. Choosing $j$ adaptively would require to restart the method from scratch whenever $j$ is changed. In such a situation it is preferable to use a sequence of rational functions with nested poles. The method of \emph{equidistributed sequences} (EDS) \cite[Section 4]{druskin2009} allows the generation of a nested sequence of poles $\xi_i$, $i=1,2,\dots$, such that the rational functions  $\prod_{i=1}^j\tfrac{z+\xi_i}{z-\xi_i}$ have asymptotically optimal convergence rates for $Z_j([-\tilde b, -\tilde a],[\tilde a,\tilde b])$ as $j$ increases, see also \cite[Section 3.5]{massei2021}. In analogy with \eqref{eq:inv-cayl}, the nested asymptotically optimal zeros and poles for $Z_j(E,F)$ are given by $\Theta^{-1}(-\xi_i)$ and $\Theta^{-1}(\xi_i)$, respectively. In the following example and more generally in Section~\ref{sec:experiments}, we demonstrate that the shift parameters computed by means of EDS yield convergence rates very close to the ones obtained with the optimal solution of  $Z_j([-\tilde b, -\tilde a],[\tilde a,\tilde b])$.
	\begin{example}\label{ex:heat1D}
		Let us consider the following $1D$ heat equation from \cite[Section 7.1]{gander2013}:
		\begin{equation}\label{eq:heat1d}
			\begin{cases}
				u_t= u_{xx} +f(x,t),&x \in \Omega:= [0,1], t \in [0,1]\\
				u\equiv 0,&\text{on }\partial \Omega,\\
				u(x,0)=4x(1-x),&\text{at } t=0,\\
				f(x,t)=h\max\{1-|c(t)-x|/w,0\},& c(t)=\frac 12 +(\frac 12-w)\sin(2\pi t), 
			\end{cases}
		\end{equation}
		with $w=0.05$  and $h=100$.
		Discretizing in space with second order central differences and in time with the implicit Euler method
		on a uniform $n\times n_t$ grid yields a matrix equation of the form $AX+XB_1^T=F$ with $A \in \R^{n\times n}$ and $B_1 \in \R^{n_t\times n_t}$ are the matrices of centered and backward differences with respect to the $x$ and $t$ variables. For this example we set $n=n_t=1000$ and rescale the equation by multiplying both sides with $\Delta t$.
		We focus on solving the correction equation
		\begin{equation}\label{eq:example1}
			A\delta X+\delta XB_1^T= -X_0\mathbf e_n\mathbf e_1^T,
		\end{equation}
		where $X_0$ satisfies $AX_0+X_0(C_1^{(1)})^T= F$. For this purpose, we use RKSM based on one of the rational functions constructed above. The poles and zeros of the rational function determine the shifts in the rational Krylov subspaces for $A$ and $B_1$, respectively, see~\eqref{eq:defrk}. We consider the following choices:
		\begin{description}
			\item[\texttt{2DISCS}] single shifts $p^*$ and $q^*$ resulting from the Zolotarev problems with 2  discs, see~\eqref{eq:singleshift};
			\item[\texttt{ZOL-DI}] optimal  shifts resulting  from  the Zolotarev problem with a disc and an interval, see~\eqref{eq:inv-cayl};
			\item[\texttt{EDS}] nested, asymptotically  optimal  shifts resulting from the Zolotarev problem with a disc and an interval; the computation of equidistributed shifts is described in \cite[Section 3.5]{massei2021}. 
			\item[\texttt{EK}] alternating shifts $0$ and $\infty$, corresponding to the \emph{extended Krylov method} \cite{simoncini2007}.
		\end{description}
		The numerical ranges of the matrices $A,B_1$, which determine the sets $E,F$, are known analytically: $\mathcal W(B_1)=\{z:\ |z+1|\leq \cos(\pi/(n_t+1))\}$ and $\mathcal W(A)=[\lambda_{\min},\lambda_{\max}]$ with $\lambda_{\min}=(2 - 2 \cos(\pi / (n+1)))\frac{(n+1)^2}{n_t}$ and $\lambda_{\max}=(2 - 2 \cos(n \pi / (n+1)))\frac{(n+1)^2}{n_t}$. We can directly use these sets 
		for generating the shift parameters in \texttt{EDS} and \texttt{ZOL-DI}, while \texttt{2DISCS} uses the (suboptimal) inclusion $\mathcal W(A)\subset \{z:\ |z-(\lambda_{\max}+\lambda_{\min})/2|\leq (\lambda_{\max}-\lambda_{\min})/2\}$.  
		
		Letting $\widetilde{\delta X}$ denote the approximate solution produced by one of the choices above, Figure~\ref{fig:krylo-conv} reports the relative error history $\norm{\widetilde{\delta X}-\delta X}_2/\norm{\delta X}_2$, where $\delta X$ is a benchmark solution computed via the \texttt{lyap} command of Matlab, as the dimension of the rational Krylov subspaces increases.
		Finally, we report the relative error associated with the truncated singular values decomposition of $\delta X$ to show the best error attainable. The results are shown in Figure~\ref{fig:krylo-conv} together with the decay bounds \eqref{eq:decay-zol-2DISCS} (\texttt{Decay 2DISCS}) and \eqref{eq:decay-zol-DI} (\texttt{Decay ZOL-DI}). The observed convergence of the residuals clearly highlights the advantage of choosing  shifts via the Zolotarev problem with a disc and an interval, either optimally (\texttt{ZOL-DI}) or asymptotically optimally (\texttt{EDS}). 
		\begin{figure}
			\centering
			\includegraphics[width=.45\textwidth]{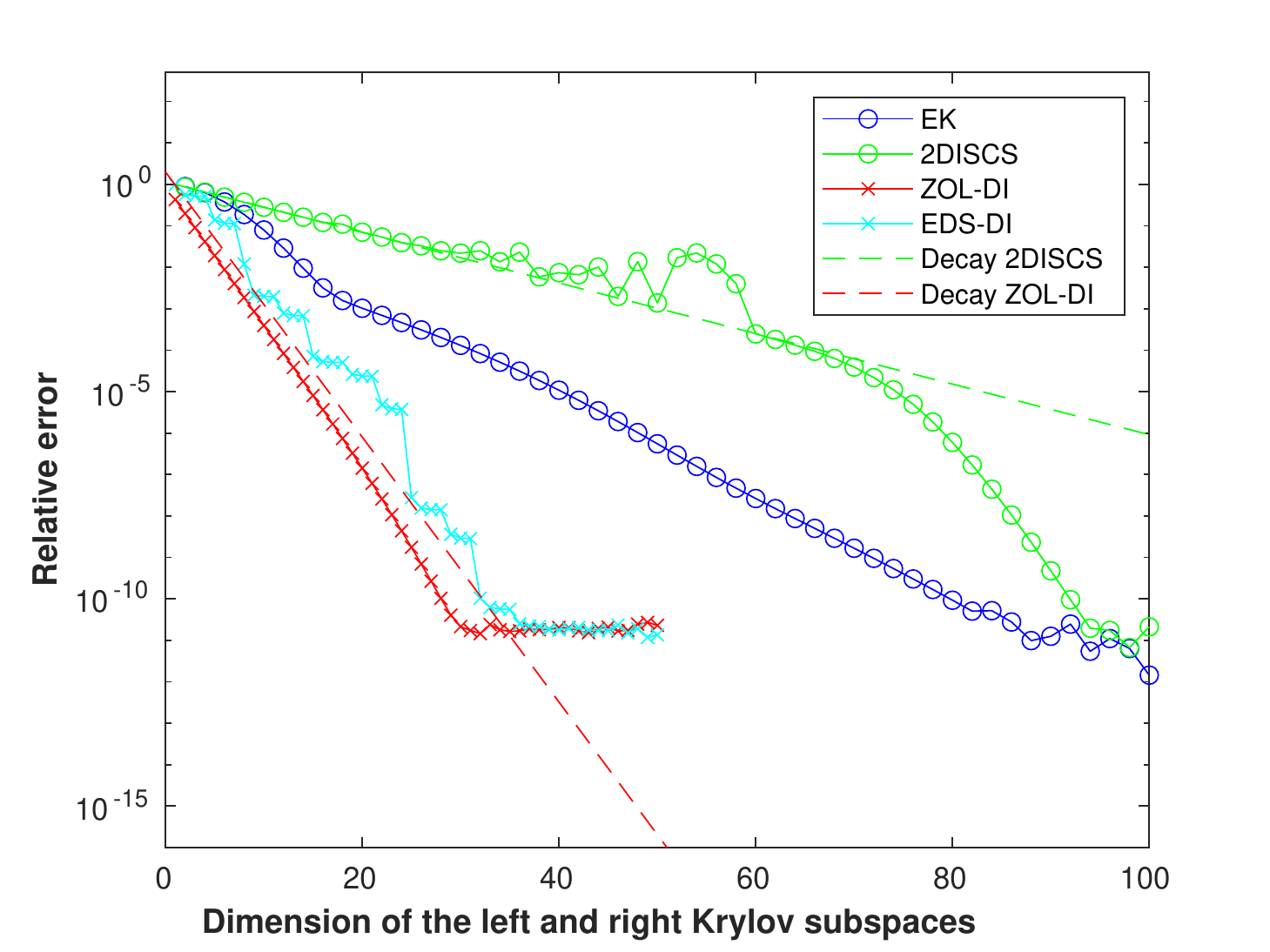}
			\caption{Relative error convergence histories of RKSM for \eqref{eq:lr-rhs}, in the case of the 1D heat equation (equation \eqref{eq:example1}), for various choices of the shift parameters.}
			\label{fig:krylo-conv}
		\end{figure}
	\end{example}
	
	\subsection{Cost analysis of Algorithm~\ref{alg:update}}\label{sec:cost}
	In order to compare with our second approach in the next section, we briefly discuss the cost of Algorithm~\ref{alg:update} as well as the potential for carrying out (embarrassingly) parallel computations.
	\upd{Let $\mathcal {C}_{\mathrm{sys}}$ denote the maximum cost of solving a linear system with a linear combination of the matrices $A$ and $M$, for a constant number of right-hand sides, respectively; note that, $\mathcal C_{\mathrm{sys}}$ is also an upper bound for the cost of solving with $A$ or $M$, individually.} The cost of a matrix-vector operation with the banded matrices $B_1,B_2$ is $\mathcal O(n_t)$, which is negligible. We may assume $k = \mathcal O(1)$, which holds for all time discretization schemes discussed in this work. We assume that RKSM converges after $\ell$ iterations to fixed accuracy and that the shift parameters are all different, as it happens, for instance, in \texttt{EDS}. If $A$ is symmetric positive definite and its condition number grows polynomially with $n$ (as, e.g., for finite difference or finite element discretizations of elliptic operators) then the bound~\eqref{eq:decay-zol-DI} predicts $\ell = \mathcal O(\log n)$. RKSM needs to solve $\mathcal O( \ell)$ linear systems at cost $\mathcal {C}_{\mathrm{sys}}$ each for generating the left factor $W$ of the approximate solution $\delta X \approx WYZ^*$. The (re)orthogonalization of the bases $W,Z$ accounts for another $\mathcal O((n+n_t)\ell^2 )$ operations. As we expect $\ell$ to be small, the cost for computing $Y$ as the solution of the compressed equation is negligible. When $M\neq I_n$, computing the low-rank factorization at line~\ref{step:UV} requires an additional constant number of linear systems solves. By adding the cost of Algorithm~\ref{alg:fast-diag} called in line 2 of Algorithm~\ref{alg:update} we obtain the overall complexity
	$$
	\mathcal O\left((n_t+\ell)\mathcal {C}_{\mathrm{sys}}+ (n+n_t)\ell^2+nn_t\log(n_t)\right).
	$$
	In terms of parallelization, note that the linear systems solves of Algorithm~\ref{alg:fast-diag} are embarrassingly parallel on up to $n_t$ cores. In this situation, the cost reduces to 
	$
	\mathcal O(\ell \mathcal {C}_{\mathrm{sys}}+ (n+n_t)\ell^2+nn_t\log(n_t))
	$
	for each of the $n_t$ cores. A further reduction can be obtained by making use of parallel implementations of RKSM~\cite{Berljafa2017} and FFT.
	
	\section{An evaluation interpolation approach}\label{sec:ev-int}
	
	As pointed out in the introduction, existing ParaDiag algorithms replace~\eqref{eq:main-eq} with 
	\begin{equation}\label{eq:all-at-once-cyclic}
		AX(C_2^{(\alpha)})^T+MX(C_1^{(\alpha)})^T=F
	\end{equation}
	for fixed $\alpha$, where $C_1^{(\alpha)},C_2^{(\alpha)}$ are $\alpha$-cyclic matrices.
	Solving equation~\eqref{eq:all-at-once-cyclic} can be carried out efficiently with Algorithm~\ref{alg:fast-diag} and
	serves as a preconditioner in a Krylov method or a stationary iteration.
	
	In this section we propose a new method that approximates the solution of \eqref{eq:main-eq} by combining solutions of the matrix equation~\eqref{eq:all-at-once-cyclic} for different values of $\alpha$. Given $\rho>0$, we let $\mathcal X_{\rho}$ denote the
	matrix-valued function that maps $z\in\mathbb C$ to the solution $X = \mathcal X_{\rho}(z)$ of~\eqref{eq:all-at-once-cyclic} for $\alpha = \rho z$. In particular $\mathcal X_{\rho}(0)$ is the desired solution of~\eqref{eq:main-eq}.
	In the following, we assume that $\mathcal X_\rho(z)$ is analytic in the disc $\{z:\ |z|<\widehat \rho\}$ for some $\widehat \rho>1$; see Example~\ref{ex:1Dheat} below for an illustration of this assumption.
	
	Our approach proceeds by interpolating $\mathcal X_\rho(z)$ at the $d$th roots of unity for some integer $d$, 
	$\omega_j=\exp(2\pi\mathbf i j/d)$ for $j=0,\dots,d-1$, with the interpolation data $\mathcal X_\rho(\omega_j)$ computed by  Algorithm~\ref{alg:fast-diag}. It is well known (see, e.g.,~\cite{trefethen2014}) that the coefficients of the interpolating polynomial $\widetilde{\mathcal X}^{(d)}(z)=\sum_{j=0}^{d-1}\widetilde X_j^{(d)} z^j$ can be obtained from the interpolation data by applying the inverse Fourier transform. In particular, we obtain
	\[
	\mathcal X_{\rho}(0) \approx \widetilde{\mathcal X}^{(d)}(0) = \widetilde X_0^{(d)} = \big(\mathcal X_\rho(\omega_0) + \cdots + \mathcal X_\rho(\omega_{d-1})\big) / d.
	\]
	The described procedure is summarized in Algorithm~\ref{alg:ev-int}. The analysis of the approximation error for such a trigonometric interpolation has a well established theory.
	The following result, retrieved by applying Theorem 12.1 from~\cite{trefethen2014} to each entry, shows that  $\widetilde X_0^{(d)}$ converges exponentially fast to the solution of~\eqref{eq:main-eq} as $d$ increases with the convergence rate depending on the radius of analyticity $\widehat \rho$. 
	\begin{theorem} \label{theorem3}
		Let $\mathcal X_\rho(z)$ be analytic on $\{z:\ |z|<\widehat \rho\}$ with $\widehat \rho >1$ and choose $R\in(1,\widehat\rho)$. Then the following holds for the entries of the approximation $\widetilde X_0^{(d)}$ returned by Algorithm~\ref{alg:ev-int}:
		$$
		|(\mathcal X_\rho(0))_{ij} -(\widetilde X_0^{(d)})_{ij}|\leq \frac{\max_{|z|=R}|(\mathcal X_\rho(z))_{ij}|}{R^d-1}, \qquad i=1,\dots n,\quad  j=1,\dots,n_t.
		$$
	\end{theorem}
	\begin{algorithm}[H] 
		\small 
		\caption{}\label{alg:ev-int}
		\begin{algorithmic}[1]
			\Procedure{ev\_int}{$A,M,B_1,B_2, F,\rho, d$}
			\For{$j=0,\dots, d-1$}
			\State $\omega_j\gets \exp(2\mathbf i j\pi / d)$
			\State $X_j \gets\textsc{fast\_diag\_solve}(A,M,B_1,B_2, F, \rho\omega_j)$ 
			\label{step:shift-solve}
			\EndFor 
			\State \Return $\widetilde X \gets (X_0 + \cdots + X_{d-1}) / d$
			\EndProcedure
		\end{algorithmic}
	\end{algorithm} 
	We remark that it is possible to implement an adaptive doubling strategy for the choice of the parameter $d$. Whenever the accuracy of $\widetilde X_0^{(d)}$  is unsatisfactory, we can consider the computation of $\widetilde X_0^{(2d)}$. Since a $d$th root of the unity is also a $2d$th root of the unity, this has the advantage that half of the evaluations of $\mathcal X_\rho(z)$ that we need in Algorithm~\ref{alg:ev-int} have already been computed at the previous iteration. As stopping criterion, one can verify that the norm
	of the difference between two consecutive approximations is smaller than a certain threshold. However, in all numerical experiments considered in this work we observe that the choices $d=2,3$ already provide sufficient accuracy and there is little need for an adaptive choice for $d$.
	\begin{example}\label{ex:1Dheat} 
		Let us consider the matrix equation $AX+XB_1^T=F$ arising from the 1D heat equation of Example~\ref{ex:heat1D} with $n=n_t$. Then $C_1^{(\rho z)} = B_1 - \rho z \mathbf e_1 \mathbf e_{n}^T$ with the matrix $B_1$ from~\eqref{eq:b1b2} multiplied with $\Delta t$.
		The eigendecomposition of $C_1^{(\rho z)}$ is given by
		$$
		C_1^{(\rho z)} =\left[\begin{smallmatrix}
			1\\
			&(\rho z)^{-\frac 1n}\\
			&&\ddots\\
			&&&(\rho z)^{\frac {1-n}n}
		\end{smallmatrix}\right]\Omega\left[\begin{smallmatrix}
			1-(\rho z)^{\frac 1n}\\
			&1-(\rho z)^{\frac 1n}\zeta\\
			&&\ddots\\
			&&&1-(\rho z)^{\frac 1n}\zeta^{n-1}
		\end{smallmatrix}\right]\Omega^*\left[\begin{smallmatrix}
			1\\
			&(\rho z)^{\frac 1n}\\
			&&\ddots\\
			&&&(\rho z)^{\frac {n-1}n}
		\end{smallmatrix}\right],
		$$
		where $\zeta=\exp(2\mathbf i \pi /n))$ and  $\Omega$ is the discrete Fourier transform.
		In particular, the eigenvalues of $C_1^{(\rho z)}$ are contained in a disc with center $1$ and radius $|\rho z|^{1/n}$. Since $A$ is symmetric positive definite, the spectra of $A$ and $-C_1^{(\rho z)}$ stay separated as long as $|z| < \rho^{-\frac1n} + \lambda_{\min}(A)$. This implies that $X_\rho(z)$ is well-defined and, as a rational function, also analytic on the open disc with radius $\rho^{\frac1n} + \lambda_{\min}(A)$ and hence the assumptions of Theorem~\ref{theorem3} can be satisfied as long as $\rho \le 1$. Assuming $|\rho z|\le 1$, the 2-norm condition number of the eigenvector matrix for $C_1^{(\rho z)}$ is bounded by $1/|\rho z|$ and, in turn,
		$\|\mathcal X_\rho(z)\|_F \le |\rho z \lambda_{\min} |^{-1} \| F \|_F$. Choosing $R = 1/\rho$ in Theorem~\ref{theorem3}, we thus obtain
		\[
		\| X_0 - \widetilde X_0^{(d)}\|_F \leq  \frac{1}{\lambda_{\min}( \rho^{-d}-1)} \|F\|_F.
		\]
		Hence, we expect to obtain very fast convergence with decreasing $\rho$ and/or increasing $d$. Note that $\rho$ cannot be chosen too small because otherwise roundoff error, due to the ill-conditioned eigenvector matrix of $X(\rho \omega_j)$, starts interfering with the interpolation error. Both statements are confirmed by numerical experiments.
		Table~\ref{tab:exp1} shows relative residual $\mathrm{Res}:=\norm{A\widetilde X_0+ \widetilde X_0B_1^T-F}_F/ \norm{F}_F$ for the matrix $\widetilde X_0$ computed by Algorithm~\ref{alg:ev-int} obtained with different values of $n,\rho,d$. The results indicate that excellent  accuracy can already achieved for $d=2$ when choosing $\rho$ properly.
		\begin{table}[ht]
			\centering  
			\caption{Relative residual of the approximate solution computed by Algorithm~\ref{alg:ev-int} with different choices of the parameters $\rho,d$ and for increasing size $n$ of the problem.}\label{tab:exp1}
			\begin{tabular}{|c|c|cccc|}
				\hline
				$n$&$\rho$& $d=1$ & $d=2$& $d=3$& $d=4$\\
				\hline
				
				\multirow{4}{*}{$500$}& $1e+00$ & $5.81e-02$ & $2.94e-06$& $1.67e-10$& $8.80e-13$ \\ 
				& $1e-02$ & $5.81e-04$ & $2.94e-10$& $3.97e-12$& $3.49e-12$ \\ 
				& $1e-04$ & $5.81e-06$ & $1.96e-10$& $1.59e-10$& $1.38e-10$ \\ 
				& $1e-06$ & $6.05e-08$ & $1.18e-08$& $9.77e-09$& $8.19e-09$ \\ 
				& $1e-08$ & $1.20e-06$ & $8.59e-07$& $7.01e-07$& $6.07e-07$ \\ 
				& $1e-10$ & $9.20e-05$ & $6.63e-05$& $5.36e-05$& $4.58e-05$ \\ 
				& $1e-12$ & $7.12e-03$ & $5.08e-03$& $4.11e-03$& $3.59e-03$ \\ 
				\hline
				\multirow{4}{*}{$1000$}& $1e+00$ & $1.59e-01$ & $6.47e-06$& $3.51e-10$& $1.61e-12$ \\ 
				& $1e-02$ & $1.59e-03$ & $6.48e-10$& $1.82e-11$& $1.59e-11$ \\ 
				& $1e-04$ & $1.59e-05$ & $1.17e-09$& $9.69e-10$& $8.31e-10$ \\ 
				& $1e-06$ & $1.93e-07$ & $7.84e-08$& $6.23e-08$& $5.53e-08$ \\ 
				& $1e-08$ & $7.74e-06$ & $5.44e-06$& $4.45e-06$& $3.83e-06$ \\ 
				& $1e-10$ & $5.78e-04$ & $4.11e-04$& $3.27e-04$& $2.84e-04$ \\ 
				& $1e-12$ & $4.50e-02$ & $3.16e-02$& $2.63e-02$& $2.19e-02$ \\  				\hline				\multirow{4}{*}{$2000$}& $1e+00$ & $1.70e-01$ & $6.70e-06$& $3.55e-10$& $7.06e-12$ \\ 
				& $1e-02$ & $1.70e-03$ & $6.75e-10$& $6.66e-11$& $5.73e-11$ \\ 
				& $1e-04$ & $1.70e-05$ & $3.79e-09$& $3.07e-09$& $2.60e-09$ \\ 
				& $1e-06$ & $3.62e-07$ & $2.25e-07$& $1.86e-07$& $1.58e-07$ \\ 
				& $1e-08$ & $2.15e-05$ & $1.52e-05$& $1.25e-05$& $1.08e-05$ \\ 
				& $1e-10$ & $1.58e-03$ & $1.12e-03$& $9.11e-04$& $7.96e-04$ \\ 
				& $1e-12$ & $1.25e-01$ & $9.00e-02$& $7.09e-02$& $6.38e-02$ \\   	\hline			\multirow{4}{*}{$4000$}& $1e+00$ & $2.14e-01$ & $9.97e-06$& $5.23e-10$& $3.96e-11$ \\ 
				& $1e-02$ & $2.14e-03$ & $1.04e-09$& $2.43e-10$& $2.12e-10$ \\ 
				& $1e-04$ & $2.14e-05$ & $1.46e-08$& $1.18e-08$& $1.03e-08$ \\ 
				& $1e-06$ & $1.34e-06$ & $9.35e-07$& $7.65e-07$& $6.63e-07$ \\ 
				& $1e-08$ & $9.56e-05$ & $6.73e-05$& $5.51e-05$& $4.78e-05$ \\ 
				& $1e-10$ & $7.32e-03$ & $5.15e-03$& $4.21e-03$& $3.66e-03$ \\ 
				& $1e-12$ & $5.98e-01$ & $4.23e-01$& $3.39e-01$& $2.98e-01$ \\ 
				\hline
			\end{tabular}
		\end{table}
	\end{example}
	\subsection{Cost analysis of Algorithm~\ref{alg:ev-int}} \label{sec:costev}
	
	In this section we maintain the assumptions and notation of Section~\ref{sec:cost}. In particular, $\mathcal {C}_{\mathrm{sys}}$ denotes the cost of solving linear systems with a linear combination of the matrices $A$ and $M$.
	The asymptotic complexity of Algorithm~\ref{alg:ev-int} is dominated by the $d$ calls to Algorithm~\ref{alg:fast-diag}, that is,
	$$
	\mathcal O(dn_t\mathcal {C}_{\mathrm{sys}} + dnn_t\log(n_t)).
	$$ 
	As the calls to Algorithm~\ref{alg:fast-diag} are embarrassingly parallel, the cost of Algorithm~\ref{alg:ev-int} reduces to 
	$$
	\mathcal O(\mathcal {C}_{\mathrm{sys}}+ nn_t\log(n_t)).
	$$
	for each core if $dn_t$ cores are available. Compared with the discussion in Section~\ref{sec:cost}, this indicates that Algorithm~\ref{alg:fast-diag} is better suited than Algorithm~\ref{alg:update} in a massively parallel environment.
	
	\section{All-at-once Runge--Kutta formulation}\label{sec:runge-kutta} 
	
	In this section, we explain how the approaches presented in this work extend when the time discretization is performed by implicit Runge--Kutta methods. For this purpose, let us consider the differential problem \eqref{eq:ode} and assume, without loss of generality, that $M=I$. The discretization by an Runge--Kutta method with $s$ stages yields
	$$
	\begin{cases}
		\mathbf x_{j+1}=\mathbf x_j+\Delta t\sum_{i=1}^sb_i \mathbf k_i^{(j)},&j=0,\dots n_t-1,\\
		\mathbf k_i^{(j)}= A\mathbf x_j+\Delta t\sum_{h=1}^sg_{ih}A \mathbf k_h^{(j)} + \mathbf f(t_j+c_i\Delta t),&i=1,\dots,s,
	\end{cases}
	$$
	where the coefficients $b_i$ and $c_i$ are the nodes and the weights of the Butcher tableau. By considering the Runge--Kutta matrix $G=(g_{ih})\in\mathbb R^{s\times s}$  and setting $K^{(j)}:=[\mathbf k^{(j)}_1|\dots| \mathbf k^{(j)}_s]\in\mathbb R^{n\times s}$, we can rewrite the second equation as
	\begin{align}\label{eq:rk1}
		\mathbf k_i^{(j)}&= A\mathbf x_j+\Delta tAK^{(j)}(G\mathbf e_i)^T + \mathbf  f(t_j+c_i\Delta t)\nonumber\\
		&\Rightarrow K^{(j)}= A\mathbf x_j\mathbf e^T+\Delta tAK^{(j)}G^T + F^{(j)},
	\end{align}
	where $\mathbf e\in\mathbb R^s$ is the vector of all ones, $\mathbf e_i\in\mathbb R^s$ is the $i$th unit vector, and $F^{(j)}:=[f(t_j+c_1\Delta t)|\dots|f(t_j+c_s\Delta t)]\in\mathbb R^{n\times s}$. Vectorizing \eqref{eq:rk1} leads to 
	$$
	\underbrace{(I-\Delta t G\otimes A)}_{H}\mathrm{vec}(K^{(j)})-(\mathbf e\otimes A)\mathbf x_j= \mathrm{vec}(F^{(j)}).
	$$ 
	Therefore, we get an all-at-once linear system $\mathcal A\mathbf x=\mathbf f$ of the form 
	\begin{equation}\label{eq:runge-kutta-lin-sys}
		\begin{bmatrix}
			H&\\
			-\mathbf b\otimes I&I\\
			& -\mathbf e\otimes A&H&\\
			&-I&-\mathbf b\otimes I&I
			\\
			&&&\ddots&\ddots&\\
			&&&&-\mathbf e\otimes A&H&\\
			&&&&-I&-\mathbf b\otimes I&I
		\end{bmatrix}
		\left[\begin{smallmatrix}
			\mathbf k_1^{(0)}\\ \vdots\\ 	\mathbf k_s^{(0)}\\ \mathbf x_1\\
			\mathbf k_1^{(1)}\\ \vdots\\	\mathbf k_s^{(1)}\\ \mathbf x_2\\
			\vdots\\ \mathbf x_{{n_t}}
		\end{smallmatrix}\right]
		=
		\left[\begin{smallmatrix}
			f(t_0+c_1\Delta t)+A\mathbf x_0\\ \vdots\\ 	f(t_0+c_s\Delta t)+A\mathbf x_0\\ \mathbf x_0\\
			f(t_1+c_1\Delta t)\\ \vdots\\	f(t_1+c_s\Delta t)\\ 0\\ f(t_2+c_1\Delta t)\\
			\vdots\\
			f(t_{n_t-1}+c_s\Delta t)\\ 0\\
		\end{smallmatrix}\right]
	\end{equation}
	where $\mathbf b:=\Delta t[b_1,\dots,b_s]$ is a row vector. Note that we can rewrite the  $2\times 2$ block matrix on the main diagonal as
	$$
	\begin{bmatrix}
		H&\\
		-\mathbf b\otimes I&I
	\end{bmatrix}=
	\begin{bmatrix}
		-\Delta t G& \phantom{e}\\
		&
	\end{bmatrix} \otimes A + \begin{bmatrix}
		I&\\
		-\mathbf b&1
	\end{bmatrix}\otimes I\in\mathbb R^{n(s+1)\times n(s+1)}.$$
	Finally, by letting $\widehat M=\left[\begin{smallmatrix}
		\phantom{0}&e\otimes A\\ &I
	\end{smallmatrix}\right]\in\mathbb R^{n(s+1)\times n(s+1)}$, we can write the coefficient matrix of the linear system \eqref{eq:runge-kutta-lin-sys} as $\mathcal A =I\otimes\widehat A + B_1\otimes \widehat M$ where
	\begin{equation}\label{eq:Acoef}
		B_1=\begin{bmatrix}
			1\\
			-1&\ddots\\
			&\ddots&\ddots\\
			&&-1&1
		\end{bmatrix},\qquad \widehat A=\begin{bmatrix}
			-\Delta t G& -\mathbf e\\
			&
		\end{bmatrix} \otimes A + \begin{bmatrix}
			I&\\
			-\mathbf b&\phantom{1}
		\end{bmatrix}\otimes I.
	\end{equation}
	This implies that~\eqref{eq:runge-kutta-lin-sys} is equivalent to a generalized Sylvester equation.
	In view of ParaDiag techniques, it is natural to consider the matrix $\mathcal A^{(\alpha)}$, obtained by replacing  $B_1$ in \eqref{eq:Acoef} with an $\alpha$-cyclic matrix $C_1^{(\alpha)}$. This leads to the following three extensions of the ParaDiag framework to Runge--Kutta methods. 
	\begin{paragraph}{Preconditioned Krylov method.}
		Employ $\mathcal A^{(\alpha)}$ as a preconditioner for GMRES.
	\end{paragraph}
	\begin{paragraph}{Evaluation interpolation approach.}
		Employ a solver for linear systems with matrices $\mathcal A^{(\rho\omega_j)}$ as building block for Algorithm~\ref{alg:ev-int}.
	\end{paragraph}
	\begin{paragraph}{Low-rank updates.}
		Solve the linear system $\mathcal A^{(1)}\mathbf x=\mathbf f$ and compute the solution of the corresponding update equation, that has right-hand side of rank $1$, with RKSM.
	\end{paragraph}
	\vspace{.2cm}
	
	\noindent In the next sections we describe the details of the solver for linear systems with the matrix $\mathcal A^{(\alpha)}$ and of the update equation that has to be addressed with RKSM.
	\subsection{Solving linear systems with $\mathcal A^{(\alpha)}$}
	After the diagonalization of the $\alpha$-cyclic matrix $C_1^{(\alpha)}$, solving linear systems with the matrix $\mathcal A^{(\alpha)}$ requires to solve --- possibly in parallel --- $n_t$ linear systems with a coefficient matrix of the form $$\lambda_j\widehat M+\widehat A=\begin{bmatrix}
		-\Delta t G& (\lambda_j-1)\mathbf e\\
		&
	\end{bmatrix} \otimes A + \begin{bmatrix}
		I&\\
		-\mathbf b&\lambda_j
	\end{bmatrix}\otimes I.$$ Each of this linear systems corresponds to solving a (generalized) Sylvester equation of the form
	\begin{equation}\label{eq:gen-sylv}
		AXN_1+XN_2= F,
	\end{equation}
	with $N_i\in\mathbb C^{(s+1)\times (s+1)}$, $i=1,2$. Because $s$, the number of Runge--Kutta stages, is usually small we can solve \eqref{eq:gen-sylv}  efficiently by computing generalized Schur decomposition~\cite{Golub2013} of $(N_1,N_2)$ via the QZ algorithm~\cite{Kag2006}. This yields unitary matrices $Q,Z\in\mathbb C^{(s+1)\times (s+1)}$ such that 
	$$
	N_1=QT_1Z^*, \qquad  N_2=QT_2Z^*,
	$$
	where $T_1,T_2\in\mathbb C^{(s+1)\times (s+1)}$ are upper triangular. The correspondingly transformed equation~\eqref{eq:gen-sylv} takes the form
	$$
	AYT_1+YT_2= \widetilde F, \qquad \widetilde F=FZ, \qquad Y=XQ,
	$$
	which can be solved by forward substitution by rewriting the equation in block-wise  form:
	\begin{align*}
		A&\begin{bmatrix}Y_1&Y_2\end{bmatrix}\begin{bmatrix}T_1^{(11)}& T_1^{(12)}\\ & T_1^{(22)}\end{bmatrix}+\begin{bmatrix}Y_1&Y_2\end{bmatrix}\begin{bmatrix}T_2^{(11)}& T_2^{(12)}\\ & T_2^{(22)}\end{bmatrix}=\begin{bmatrix}\widetilde F_1&\widetilde F_2\end{bmatrix}\\ &\Longrightarrow\begin{cases}
			AY_1T_1^{(11)}+ Y_1T_2^{(11)}=\widetilde F_1\\
			AY_2T_1^{(22)}+ Y_2T_2^{(22)}=\widetilde F_2-AY_1T_1^{(12)}-Y_1T_2^{(12)}
		\end{cases}.
	\end{align*}
	Assuming  that $T_1^{(11)}, T_2^{(11)}$ are scalar quantities we can first retrieve $Y_1\in \mathbb C^n$ by solving a linear system with $(A+T_2^{(11}) I)$ and continuing to compute $Y_2\in \mathbb C^{n\times s}$ via recursion. The solution $X$ of~\eqref{eq:gen-sylv} is obtained as $YQ^*$.
	
	\subsection{Solving the update equation}
	In the case of Runge--Kutta methods, the low-rank update approach from Section~\ref{sec:low-rank} requires us to solve the generalized matrix equation
	\begin{equation}\label{eq:lowrank-sylv}
		\widehat A\delta X +\widehat M\delta XB_1^T= \widehat MX_0\mathbf e_{n_t}\mathbf e_1^T,
	\end{equation}
	where $X_0\in\mathbb C^{n(s+1)\times n_t}$ is the (reshaped) solution of the linear system $\mathcal A^{(1)}\mathbf x=\mathbf f$  described in the previous section. Equation~\eqref{eq:lowrank-sylv} has a right-hand side of rank $1$ but it is not immediate to recast it as a  Sylvester equation of the form \eqref{eq:lr-rhs} because the coefficient $\widehat M$ is not invertible. Given a shift parameter $\sigma\in\mathbb C$, we add and subtract the quantity $\sigma\widehat A\delta XB_1^T$ to \eqref{eq:lowrank-sylv}, obtaining
	\[
	\widehat A\delta X(I-\sigma B_1^T) +(\widehat M+\sigma \widehat A)\delta XB_1^T=\widehat MX_0\mathbf e_{n_t}\mathbf e_1^T.
	\]
	Hence, if $\sigma$ is chosen such that both $\widehat M+\sigma \widehat A$ and $I-\sigma B_1^T$ are invertible, we can write an equation of the form $\widetilde A \delta X+\delta X\widetilde B^T = UV^*$ equivalent to \eqref{eq:lowrank-sylv} by setting:
	\begin{align*}
		\widetilde A&= (\widehat M+\sigma \widehat A)^{-1}\widehat A,&
		\widetilde B&=(I-\sigma B_1)^{-1}B_1,\\
		\mathbf u &=(\widehat M+\sigma \widehat A)^{-1}\widehat MX_0\mathbf e_{n_t},&
		\mathbf v &= (I-\sigma B_1)^{-1}\mathbf e_1.
	\end{align*}
	This suggests to use the matrices $\widetilde A,\widetilde B$ and the vectors $\mathbf u, \mathbf v$ to generate Krylov subspaces for the equation \eqref{eq:lowrank-sylv}. Once again, there is no need to form explicitly $\widetilde A,\widetilde B$ for performing matrix-vector products and solving shifted linear systems; for these operations we rely on the relations:
	\begin{align*}
		\left[(\widehat M+\sigma \widehat A)^{-1}\widehat A-zI\right]^{-1}\mathbf w&= \left[(1-z\sigma)\widehat A-z\widehat M\right]^{-1}(\widehat M+\sigma \widehat A)\mathbf w,\\
		\left[(I-\sigma B_1)^{-1}B_1-zI\right]^{-1}\mathbf w &=\left[(1+z\sigma)B_1-zI\right]^{-1}(I-\sigma B_1)\mathbf w, 
	\end{align*}
	that only need to manipulate the matrices $\widehat A, \widehat M, B_1,I-\sigma B_1$, which are sparse  when $A$ is sparse.
	
	In order to study the low-rank property of the solution to~\eqref{eq:lowrank-sylv} and to provide suitable a priori choices for the shift parameters in RKSM, one would need to study the spectral properties of $\widetilde A,\widetilde B$. The matrix  $\widetilde B$ is Toeplitz triangular and we can explicitly compute its entries:
	$$
	\widetilde B=(1-\sigma)^{-1}
	\begin{bmatrix}
		1\\
		\frac{1}{\sigma-1}&\ddots\\
		\frac{\sigma}{(\sigma-1)^2}&\ddots&\ddots\\
		\vdots&\ddots&\ddots&\ddots\\
		\frac{\sigma^{n_t-2}}{(\sigma-1)^{n_t-1}}&\dots&\frac{\sigma}{(\sigma-1)^2}&\frac{1}{\sigma-1}&1
	\end{bmatrix}
	$$
	In particular, it has the only eigenvalue $(1-\sigma)^{-1}$ with multiplicity $n_t$. However, the matrix $\widetilde A$ appears hard to analyze in general, even assuming strong properties such as $A$ symmetric positive definite. For this reason, we postpone the analysis of \eqref{eq:lowrank-sylv} to future work and we propose the use of the extended Krylov method (see Example~\ref{ex:1Dheat}) or adaptive choices of the shift parameters \cite{druskin2011} when solving the update equation with RKSM.  
	
	\section{Numerical results}\label{sec:experiments}
	
	In this section, we compare the performances of Algorithm~\ref{alg:update} and Algorithm~\ref{alg:ev-int} with those of the preconditioned GMRES method considered  in \cite{mcdonald2018,gander2020paradiag}, on different examples.  More specifically, we consider the following list of methods:
	\begin{itemize}
		\item \texttt{PGMRES}: GMRES preconditioned by the solver of the matrix equation with a $\alpha$-circulant coefficient. If not stated otherwise, $\alpha=1$ and the threshold for stopping GMRES is set to $10^{-8}$,
		\item \texttt{Ev-Int}: Algorithm~\ref{alg:ev-int}. If not stated otherwise, we choose $d=2$ and $\rho=5\cdot 10^{-4}$.
		\item \texttt{2DISCS}, \texttt{ZOL-DI(4)}, \texttt{EDS}, \texttt{EK}: Algorithm~\ref{alg:update} with different choices of poles; see Example~\ref{ex:heat1D} for details. Note that, the method \texttt{ZOL-DI} is not competitive, in terms of computational time, as it requires to recompute the bases of the Krylov subspaces from scratch at every iteration. For this reason we consider its variant \texttt{ZOL-DI(4)}, that consists in cyclically repeating the $4$ shifts corresponding to the Zolotarev problem with rational functions of degree $(4,4)$. We have chosen $4$ shifts because preliminary experiments indicated that this choice offers good performance across a broader range of examples.
	\end{itemize}
	RKSM called at line~\ref{step:rksm} of Algorithm~\ref{alg:update} is stopped when the relative residual of the update equation~\eqref{eq:lr-rhs} is below the threshold $10^{-8}$.  The methods \texttt{EK},  \texttt{2DISCS}, and \texttt{ZOL-DI(4)} rely on precomputed factorizations of the, possibly shifted, matrices $\widetilde A$ and $\widetilde B$ when generating the bases of the rational Krylov subspaces.  
	If the matrices $A$ and $M$ in \eqref{eq:main-eq} are sparse, a sparse direct solver is exploited whenever the solution of a linear system involving a linear combination of $A$ and $M$ is required; e.g., at line~\ref{step:shift-solve} of Algorithm~\ref{alg:ev-int}.  
	
	For each method we report the computational time required and the relative error $\mathrm{Res}:= \norm{A\widehat XB_2^T+M\widehat XB_1^T-F}_F/\norm{F}_F$ in order to assess the accuracy of the correspondent approximate solution $\widehat X$. \upd{The best timings of each case study are highlighted in bold.} For the various implementations of Algorithm~\ref{alg:update} we also report the percentage of the CPU time that is spent for solving \eqref{eq:lr-rhs}, labeled with $\mathrm{T_{sylv}}$, the dimension (Dim) of the rational Krylov subspace generated and the rank (rk) of the approximate solution to \eqref{eq:lr-rhs}, after recompression.
	
	\upd{Apart from Section~\ref{sec:parallel-exp},} all experiments have been performed on a laptop with a dual-core Intel Core i7-7500U 2.70 GHz CPU, 256KB of level 2 cache, and 16 GB of RAM. The algorithms
	are implemented in MATLAB and tested under MATLAB2020a, with MKL BLAS
	version 2019.0.3 utilizing both cores.
	
	\subsection{1D heat equation}
	
	We start by testing all the algorithms detailed at the beginning of this section on the matrix equation $AX+XB_1^T=F$ from Example~\ref{ex:heat1D} for $n=1089, 4225$ and $n_t=2^8,\dots, 2^{11}$. 
	
	The results are reported in Table~\ref{tab:heat1D}. As expected \texttt{ZOL-DI(4)} and \texttt{EDS} generate significantly lower dimensional Krylov subspaces with respect to those generated by \texttt{EK} and \texttt{2DISCS}, while achieving comparable accuracies. Looking at the computational times, we see that, for all methods, we have an almost linear dependence on the size $n$ and on the number of time steps $n_t$.  The methods with the most favorable timings are \texttt{Ev-Int} and \texttt{ZOL-DI(4)}. The former is preferable when $\mathrm{T_{sylv}}>50 \%$, which  happens for the smaller instances of the equation. For all four implementations of Algorithm~\ref{alg:update}, the fraction of the time spent for solving  the low-rank Sylvester equation decreases as the size of the problem increases. Since the advantage of using \texttt{ZOL-DI(4)} is inversely proportional to $\mathrm{T_{sylv}}$ we get the best scenario for \texttt{ZOL-DI(4)} when $n=4225$ and $n_t=2048$ where we get about a $2\times$ speedup over \texttt{Ev-Int} and a $5\times$ speedup over \texttt{PGMRES}.

	\subsection{Convection diffusion}
	
	This example is concerned with the convection diffusion equation from \cite[Section 6.2]{mcdonald2018}, which models the temperature distribution in a cavity having an external wall maintained at a constant (hot) temperature and a wind that determines a recirculating ﬂow. The problem reads as follows:
	$$
	\begin{cases}
		u_t-\epsilon\Delta u + \mathbf w\cdot \nabla u=0,&S\times (0,1),\ S = (-1,1)^2,\\
		u(x, y, 0)=\chi_{\{x=1\}},\\
		u(x, y, t)=\chi_{\{x=1\}},& \partial S\times (0,1),
	\end{cases}
	$$
	where $\epsilon = 0.005$, $\chi_{\{x=1\}}$ is the indicator function of the set $\{x=1\}$, and $\mathbf w = (2y(1-x^2), -2x(1-y^2))$. The discretization is done via Q1 finite elements over the domain $S$ and backward Euler time-stepping, together with streamline-upwind Petrov--Galerkin (SUPG) stabilization. The matrix equation corresponding to the all-at-once linear system takes the form $AX+MXB_1^T=F$, where both $A$ and $M$ are sparse. To estimate the numerical range of $M^{-1}A$, we compute the largest and smallest real parts, $r_{\max},r_{\min}$,  of its eigenvalues, by means of the Matlab command \texttt{eigs}. We remark that, $r_{\max},r_{\min}>0$ and we denote $c:=(r_{\max}+r_{\min})/2$. Then, \texttt{2DISCS} employs $\{z:\ |z-c|\leq (r_{\max}-r_{\min})/2 \}$ as estimate for $\mathcal W(M^{-1}A)$ while \texttt{ZOL-DI(4)} and \texttt{EDS} makes use of the real interval $[r_{\min},r_{\max}]$. The performances of the numerical methods are reported in Table~\ref{tab:2dconv-diff}. Although we used the slightly higher tolerance $10^{-6}$ to stop \texttt{PGMRES}, the latter needs 17 iterations to converge and this makes the other approaches significantly faster. Among the methods based on low-rank updates, \texttt{ZOL-DI(4)} and \texttt{EDS} are those that generate the Krylov subspaces of smallest dimensions, while maintaining a good accuracy. This feature makes \texttt{EDS} faster than \texttt{EK} although the latter has a cheaper iteration cost. Similar comments to the ones made in the previous example apply to the comparison between $\texttt{ZOL-DI(4)}$ and $\texttt{Ev-Int.}$
	
	\subsection{Fractional space diffusion}\label{sec:frac-space}
	
	We consider the fractional diffusion example from \cite[Section 5]{wu2018}:
	$$
	\begin{cases}
		u_t= a_1\frac{\partial^{\gamma_1} u^+}{\partial x^{\gamma_1}}+b_1\frac{\partial^{\gamma_1} u^-}{\partial x^{\gamma_1}}+a_2\frac{\partial^{\gamma_2} u^+}{\partial y^{\gamma_2}}+b_2\frac{\partial^{\gamma_2} u^-}{\partial y^{\gamma_2}}+f(x,y,t), & S \times (0,5),\ S = (0,1)^2,\\
		u(x,y, 0)= 0, & (x,y)\in S,\\
		u(x,y,t)=0,& (x,y,t)\in\partial S\times (0,5),
	\end{cases}
	$$
	where $\frac{\partial^{\gamma_j}u^\pm}{\partial x}$ indicates the Riemann-Liouville right- and left- looking derivatives; the diffusion coefficients are set as $a_1 = 1, a_2 = 0.5, 
	b_1 = 0.2, b_2 = 1$, the differential orders are $\gamma_1=1.75$ and $\gamma_2=1.5$, and the  source term is $f(x,y,t)= 10\sin(3xyt)$. For space discretization, the second-order weighted and shifted Gr\"unwald
	difference formula is employed; uniform backward Euler time stepping is applied. The all-at-once linear system matrix takes the form $I\otimes A+B_1\otimes I$ where the matrix $A$ has the structure
	$$
	A = I \otimes T_1+T_2\otimes I
	$$
	with $T_1,T_2$ non sparse Toeplitz matrices; see \cite{wu2018} for a detailed description. In particular, matrix-vector operations with the matrix $A$ are efficiently performed by relying on FFT and matrix equation techniques. Moreover, the minimum and maximum real part of the eigenvalues of $A$ are easily obtained once estimates of the same quantities for the matrices $T_1,T_2$ are computed. When running $\texttt{2DISCS}, \texttt{ZOL-DI(4)}, \texttt{EDS}$, we approximate $\mathcal W(A)$ as we did in the previous example.
	
\begin{figure}
		\centering
		\includegraphics[width=.4\textwidth]{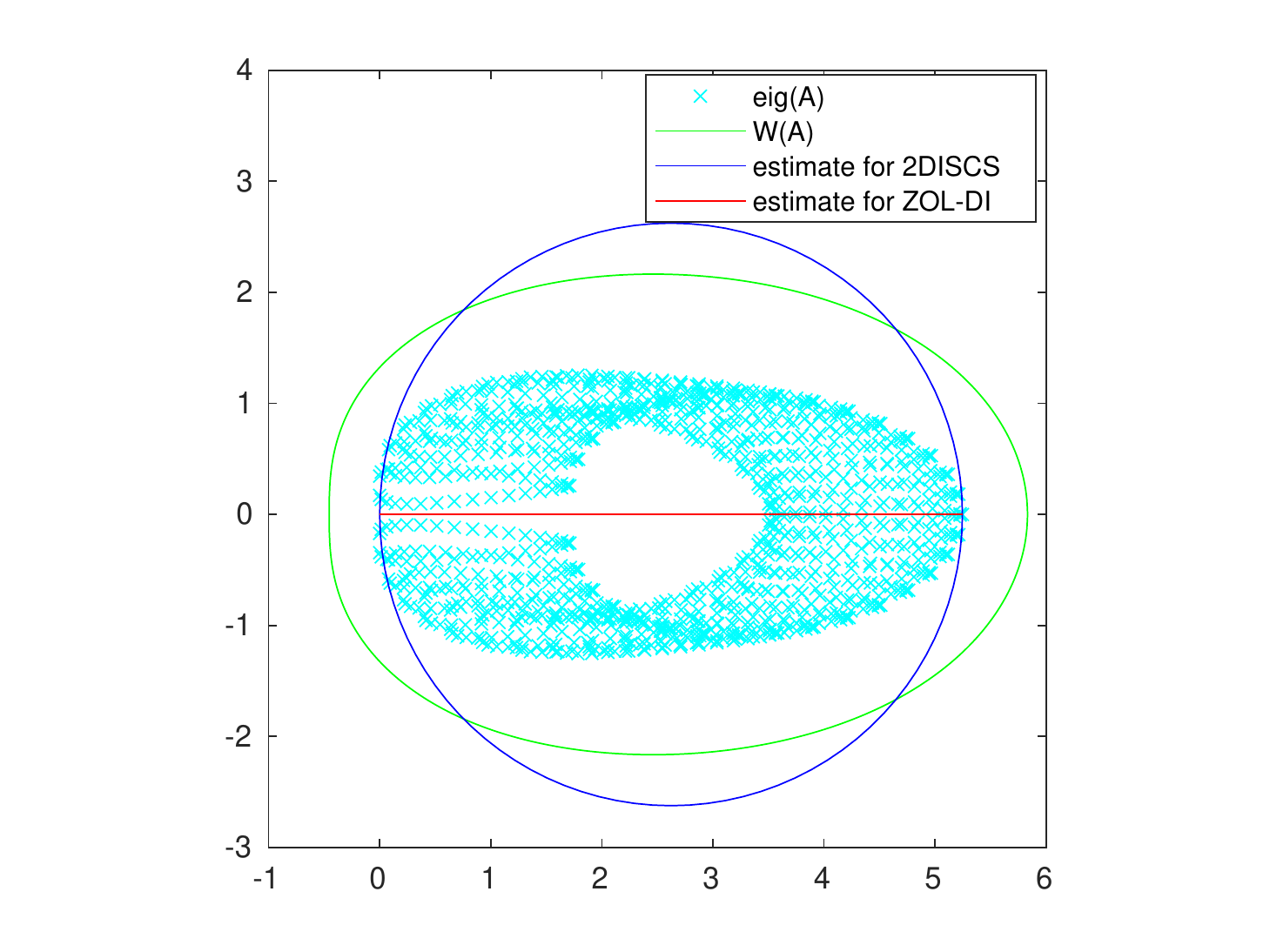}
		\caption{Fractional diffusion example with $n=1089$; eigenvalues of the matrix $A$ (cyan crosses), numerical range of the matrix $A$ (disc enclosed by the green line) and its approximations used for \texttt{2DISCS} (disc enclosed by the blue line), and for \texttt{ZOL-DI(4)} (red line segment).}
		\label{fig:numericalrange}
	\end{figure}

	The results reported in Table~\ref{tab:fract-diff} show that the methods based on low-rank updates struggle on this example because they have to generate somewhat large subspaces in order to achieve the desired accuracy. This is caused by the fact that the solution of the update equation \eqref{eq:update-eq} has a slightly higher rank with respect to the previous examples, and the approximations of $\mathcal W(A)$ are not tight, see Figure~\ref{fig:numericalrange}. Nevertheless, \texttt{ZOL-DI(4)} achieves the best timings for $n=4225$. \texttt{Ev-Int} maintains the usual behavior and turns out to be the preferable method for $n=1089$. 
	\begin{landscape}
		\begin{table}[ht]
			\centering
			\caption{1D heat equation; performances of the various methods.}		\label{tab:heat1D}
			\resizebox{1.4\textheight}{!}{  
				\begin{tabular}{|c|c|ccc|cc|ccccc|ccccc|ccccc|ccccc|}
					\hline
					&&\multicolumn{3}{c}{\texttt{PGMRES}}&\multicolumn{2}{|c|}{\texttt{Ev-Int}}&\multicolumn{5}{|c|}{\texttt{EK}}&\multicolumn{5}{|c|}{\texttt{2DISCS}}&\multicolumn{5}{|c|}{\texttt{EDS}}&\multicolumn{5}{|c|}{\texttt{ZOL-DI(4)}}\\
					\hline
					$n$&$n_t$& Time &Res&It&Time&Res&Time&Res&rk&$T_{sylv}$ (\%)&Dim&Time&Res&rk&$T_{sylv}$ (\%)&Dim&Time&Res&rk&$T_{sylv}$ (\%)&Dim&Time&Res&rk&$T_{sylv}$ (\%)&Dim\\
					\hline
					
					\multirow{4}{*}{$1089$}& $256$ & $0.25$ & $7.2e-10$ & $3$ & $\mathbf{0.10}$ & $2.4e-09$ & $0.39$ & $1.4e-08$ & $12$ & $88.3$ & $47$ & $0.30$ & $8.9e-09$ & $12$ & $85.1$ & $49$ & $0.24$ & $2.6e-10$ & $12$ & $81.3$ & $25$ & $0.15$ & $4.7e-11$ & $12$ & $69.7$ & $33$ \\ 
					& $512$ & $0.47$ & $2.1e-09$ & $3$ & $\mathbf{0.21}$ & $2.3e-09$ & $0.48$ & $1.0e-08$ & $13$ & $78.0$ & $59$ & $0.57$ & $5.1e-09$ & $13$ & $81.6$ & $67$ & $0.27$ & $7.8e-09$ & $13$ & $61.4$ & $25$ & $0.22$ & $7.8e-11$ & $13$ & $51.9$ & $33$ \\ 
					& $1024$ & $0.98$ & $1.6e-09$ & $3$ & $0.40$ & $2.4e-09$ & $0.70$ & $2.5e-08$ & $14$ & $75.4$ & $71$ & $1.00$ & $2.1e-08$ & $14$ & $82.7$ & $89$ & $0.37$ & $1.6e-08$ & $14$ & $53.7$ & $29$ & $\mathbf{0.33}$ & $3.0e-09$ & $14$ & $47.4$ & $33$ \\ 
					& $2048$ & $2.02$ & $3.5e-09$ & $3$ & $0.82$ & $2.3e-09$ & $1.14$ & $3.3e-08$ & $15$ & $68.3$ & $89$ & $1.91$ & $2.8e-08$ & $15$ & $80.9$ & $121$ & $0.61$ & $9.7e-09$ & $15$ & $39.9$ & $32$ & $\mathbf{0.51}$ & $5.1e-09$ & $15$ & $28.6$ & $37$ \\  
					\hline
					\multirow{4}{*}{$4225$}& $256$ & $1.08$ & $2.9e-08$ & $3$ & $0.40$ & $3.4e-08$ & $0.54$ & $1.4e-08$ & $13$ & $70.2$ & $47$ & $0.63$ & $9.0e-09$ & $13$ & $74.4$ & $49$ & $0.41$ & $7.2e-10$ & $13$ & $60.2$ & $25$ & $\mathbf{0.33}$ & $6.8e-10$ & $13$ & $50.9$ & $33$ \\ 
					& $512$ & $2.27$ & $3.0e-08$ & $3$ & $0.75$ & $3.4e-08$ & $0.86$ & $1.0e-08$ & $14$ & $62.4$ & $59$ & $1.08$ & $5.2e-09$ & $14$ & $69.5$ & $67$ & $0.59$ & $7.9e-09$ & $14$ & $43.4$ & $25$ & $\mathbf{0.51}$ & $7.1e-10$ & $14$ & $34.7$ & $33$ \\ 
					& $1024$ & $4.56$ & $1.3e-07$ & $3$ & $1.69$ & $3.7e-08$ & $1.66$ & $2.5e-08$ & $15$ & $52.0$ & $71$ & $2.21$ & $2.1e-08$ & $15$ & $63.6$ & $89$ & $1.25$ & $1.6e-08$ & $15$ & $36.3$ & $29$ & $\mathbf{0.96}$ & $3.1e-09$ & $15$ & $17.4$ & $33$ \\ 
					& $2048$ & $9.72$ & $3.1e-07$ & $3$ & $3.28$ & $3.6e-08$ & $2.90$ & $3.3e-08$ & $17$ & $48.3$ & $89$ & $4.07$ & $2.9e-08$ & $17$ & $62.6$ & $121$ & $1.85$ & $9.7e-09$ & $17$ & $18.4$ & $32$ & $\mathbf{1.79}$ & $5.0e-09$ & $17$ & $15.7$ & $37$ \\ 
					\hline
				\end{tabular}
			}
		\end{table}
		\vspace{1cm}
		
		\begin{table}[h!]
			\centering
			\caption{Convection diffusion; performances of the various methods.}
			\label{tab:2dconv-diff}
			\resizebox{1.4\textheight}{!}{  
				\begin{tabular}{|c|c|ccc|cc|ccccc|ccccc|ccccc|ccccc|}
					\hline
					&&\multicolumn{3}{c}{\texttt{PGMRES}}&\multicolumn{2}{|c|}{\texttt{Ev-Int}}&\multicolumn{5}{|c|}{\texttt{EK}}&\multicolumn{5}{|c|}{\texttt{2DISCS}}&\multicolumn{5}{|c|}{\texttt{EDS}}&\multicolumn{5}{|c|}{\texttt{ZOL-DI(4)}}\\
					\hline
					$n$&$n_t$& Time &Res&It&Time&Res&Time&Res&rk&$T_{sylv}$ (\%)&Dim&Time&Res&rk&$T_{sylv}$ (\%)&Dim&Time&Res&rk&$T_{sylv}$ (\%)&Dim&Time&Res&rk&$T_{sylv}$ (\%)&Dim\\
					\hline
					
					\multirow{4}{*}{$1089$}& $256$ & $19.66$ & $4.2e-08$ & $17$ & $2.05$ & $7.2e-10$ & $3.45$ & $2.3e-10$ & $28$ & $70.6$ & $165$ & $2.66$ & $2.2e-10$ & $28$ & $61.9$ & $123$ & $2.03$ & $1.6e-10$ & $28$ & $50.0$ & $67$ & $\mathbf{1.49}$ & $1.8e-10$ & $28$ & $32.0$ & $65$ \\ 
					& $512$ & $42.25$ & $7.0e-08$ & $17$ & $4.29$ & $1.0e-09$ & $4.90$ & $4.4e-11$ & $29$ & $70.2$ & $171$ & $4.22$ & $5.3e-11$ & $29$ & $52.0$ & $141$ & $3.13$ & $2.0e-11$ & $29$ & $41.6$ & $70$ & $\mathbf{2.72}$ & $2.3e-11$ & $29$ & $30.8$ & $69$ \\ 
					& $1024$ & $81.87$ & $1.0e-07$ & $17$ & $8.62$ & $1.5e-09$ & $6.79$ & $1.7e-10$ & $30$ & $55.7$ & $165$ & $6.79$ & $2.5e-10$ & $30$ & $48.5$ & $153$ & $5.29$ & $2.7e-10$ & $30$ & $30.8$ & $72$ & $\mathbf{4.84}$ & $1.8e-10$ & $30$ & $20.1$ & $69$ \\ 
					& $2048$ & $168.41$ & $1.4e-07$ & $17$ & $17.00$ & $2.1e-09$ & $11.20$ & $4.7e-11$ & $31$ & $55.0$ & $165$ & $11.20$ & $6.4e-11$ & $31$ & $40.9$ & $159$ & $9.57$ & $8.7e-11$ & $31$ & $25.9$ & $75$ & $\mathbf{9.11}$ & $3.6e-11$ & $31$ & $18.9$ & $77$ \\ 
					\hline
					\multirow{4}{*}{$4225$}& $256$ & $92.15$ & $1.7e-08$ & $17$ & $9.59$ & $8.4e-11$ & $15.95$ & $4.7e-10$ & $31$ & $36.9$ & $251$ & $9.90$ & $3.4e-10$ & $32$ & $37.0$ & $147$ & $8.13$ & $3.1e-10$ & $31$ & $19.2$ & $82$ & $\mathbf{6.87}$ & $3.0e-10$ & $31$ & $11.6$ & $81$ \\ 
					& $512$ & $191.68$ & $2.6e-08$ & $17$ & $18.45$ & $1.2e-10$ & $20.54$ & $7.1e-11$ & $33$ & $38.3$ & $259$ & $15.65$ & $7.2e-11$ & $33$ & $29.0$ & $171$ & $12.48$ & $6.5e-11$ & $33$ & $16.5$ & $84$ & $\mathbf{11.39}$ & $3.0e-11$ & $33$ & $10.9$ & $85$ \\ 
					& $1024$ & $367.04$ & $3.9e-08$ & $17$ & $36.18$ & $1.7e-10$ & $29.36$ & $5.4e-10$ & $34$ & $24.6$ & $259$ & $25.52$ & $6.4e-10$ & $35$ & $24.6$ & $189$ & $21.70$ & $3.6e-10$ & $35$ & $11.7$ & $92$ & $\mathbf{20.34}$ & $1.9e-10$ & $35$ & $7.2$ & $89$ \\ 
					& $2048$ & $725.21$ & $5.6e-08$ & $17$ & $72.69$ & $2.4e-10$ & $47.54$ & $1.7e-10$ & $35$ & $23.3$ & $251$ & $44.69$ & $1.1e-10$ & $35$ & $18.4$ & $197$ & $40.20$ & $1.4e-10$ & $35$ & $9.3$ & $93$ & $\mathbf{38.75}$ & $7.1e-11$ & $35$ & $5.9$ & $89$ \\ 
					\hline
				\end{tabular}
			}
		\end{table}
		\vspace{1cm}
		
		\begin{table}[h!]
			\centering
			\caption{Fractional space diffusion; performances of the various methods.}
			\label{tab:fract-diff}
			\resizebox{1.4\textheight}{!}{  
				\begin{tabular}{|c|c|ccc|cc|ccccc|ccccc|ccccc|ccccc|}
					\hline
					&&\multicolumn{3}{c}{\texttt{PGMRES}}&\multicolumn{2}{|c|}{\texttt{Ev-Int}}&\multicolumn{5}{|c|}{\texttt{EK}}&\multicolumn{5}{|c|}{\texttt{2DISCS}}&\multicolumn{5}{|c|}{\texttt{EDS}}&\multicolumn{5}{|c|}{\texttt{ZOL-DI(4)}}\\
					\hline
					$n$&$n_t$& Time &Res&It&Time&Res&Time&Res&rk&$T_{sylv}$ (\%)&Dim&Time&Res&rk&$T_{sylv}$ (\%)&Dim&Time&Res&rk&$T_{sylv}$ (\%)&Dim&Time&Res&rk&$T_{sylv}$ (\%)&Dim\\
					\hline
					
					\multirow{4}{*}{$1089$}& $256$ & $\mathbf{2.03}$ & $3.9e-12$ & $2$ & $0.92$ & $4.9e-11$ & $2.35$ & $1.3e-07$ & $26$ & $74.4$ & $113$ & $1.45$ & $3.8e-08$ & $26$ & $58.5$ & $97$ & $2.37$ & $1.3e-07$ & $26$ & $74.7$ & $88$ & $1.34$ & $8.9e-08$ & $26$ & $53.8$ & $85$ \\ 
					& $512$ & $\mathbf{3.42}$ & $1.3e-10$ & $2$ & $1.68$ & $8.1e-11$ & $4.36$ & $1.6e-07$ & $33$ & $79.4$ & $171$ & $2.41$ & $5.8e-08$ & $33$ & $62.7$ & $133$ & $3.96$ & $1.9e-07$ & $33$ & $77.3$ & $121$ & $2.14$ & $1.6e-07$ & $33$ & $57.9$ & $117$ \\ 
					& $1024$ & $\mathbf{8.74}$ & $1.2e-10$ & $3$ & $3.69$ & $1.0e-10$ & $13.15$ & $1.7e-07$ & $40$ & $86.4$ & $265$ & $4.69$ & $2.5e-07$ & $40$ & $62.0$ & $169$ & $9.24$ & $2.8e-07$ & $40$ & $80.7$ & $175$ & $4.54$ & $1.4e-07$ & $40$ & $60.7$ & $169$ \\ 
					& $2048$ & $\mathbf{21.64}$ & $1.2e-11$ & $4$ & $7.14$ & $5.8e-08$ & $29.90$ & $3.3e-07$ & $46$ & $87.8$ & $355$ & $7.59$ & $3.9e-07$ & $46$ & $52.0$ & $197$ & $16.85$ & $1.9e-07$ & $46$ & $78.4$ & $220$ & $8.61$ & $3.9e-07$ & $46$ & $57.6$ & $213$ \\
					\hline
					\multirow{4}{*}{$4225$}& $256$ & $7.30$ & $3.5e-11$ & $2$ & $3.58$ & $4.2e-10$ & $3.78$ & $1.0e-07$ & $25$ & $51.6$ & $107$ & $3.48$ & $1.2e-07$ & $25$ & $47.4$ & $97$ & $4.73$ & $2.1e-07$ & $25$ & $61.3$ & $89$ & $\mathbf{3.09}$ & $3.1e-08$ & $25$ & $40.7$ & $89$ \\ 
					& $512$ & $14.47$ & $7.2e-11$ & $2$ & $7.15$ & $5.5e-10$ & $10.01$ & $2.4e-07$ & $32$ & $63.8$ & $179$ & $6.74$ & $1.1e-07$ & $32$ & $46.1$ & $141$ & $9.42$ & $2.6e-07$ & $32$ & $61.5$ & $126$ & $\mathbf{6.23}$ & $3.0e-07$ & $32$ & $41.6$ & $121$ \\ 
					& $1024$ & $35.71$ & $1.9e-08$ & $2$ & $14.16$ & $4.3e-10$ & $22.50$ & $3.0e-07$ & $42$ & $68.7$ & $255$ & $13.86$ & $2.4e-07$ & $42$ & $49.2$ & $197$ & $16.70$ & $2.8e-07$ & $42$ & $57.9$ & $170$ & $\mathbf{11.23}$ & $2.9e-07$ & $42$ & $37.3$ & $165$ \\ 
					& $2048$ & $88.22$ & $1.0e-08$ & $3$ & $29.11$ & $1.6e-09$ & $54.76$ & $5.9e-07$ & $55$ & $71.8$ & $375$ & $26.86$ & $4.4e-07$ & $55$ & $42.4$ & $265$ & $41.61$ & $5.5e-07$ & $55$ & $62.8$ & $256$ & $\mathbf{24.59}$ & $5.1e-07$ & $55$ & $37.1$ & $237$ \\
					\hline
			\end{tabular}}
		\end{table}
	\end{landscape}
	\subsection{Fractional time diffusion}
	
	We now consider an example where Algorithm~\ref{alg:update} is not applicable due to the non locality of the time differential operator. We modify Example~\ref{ex:1Dheat} by replacing the first-order time derivative with the Caputo time derivative of order $\gamma\in(0,1)$:
	$$
	\begin{cases}
		\frac{\partial^\gamma u}{\partial t^\gamma}= u_{xx} +f(x,t),&S:= [0,1]\times [0, 1],\\
		u\equiv 0,&\text{on }\partial S,\\
		u(x,0)=4x(1-x),&\text{at } t=0,\\
		f(x,t)=h\max\{1-|c(t)-x|/w,0\},& c(t)=\frac 12 +(\frac 12-w)\sin(2\pi t), 
	\end{cases}
	$$
	and we keep the same parameters $w,h$ of the integer derivative case.  Discretizing the Caputo time derivative with the unshifted Grünwald-Letnikov formula, we get a matrix equation $AX+XB_1^T=F$ where $B_1$ is Toeplitz lower triangular and takes the form:
	$$
	B_1=\begin{bmatrix}
		g_{\gamma,0}\\
		g_{\gamma, 1}&g_{\gamma,0}\\
		\vdots&\ddots&\ddots\\
		g_{\gamma, n_t-1}&\dots&g_{\gamma, 1}&g_{\gamma, 0}
	\end{bmatrix}, \qquad \begin{cases}
		g_{\gamma,0}=1,\\
		g_{\gamma,j}= \left(1-\frac{\gamma+1}{j}\right)g_{\gamma, j-1}.
	\end{cases}
	$$
	In particular, the usual way to construct a circulant matrix from $B_1$ requires to apply a matrix $\delta B_1^{(\alpha)}$ of rank $n_t-1$. For this reason we restrict to \texttt{PGMRES} and \texttt{Ev-Int} for solving the all-at-once problem. The results reported in Table~\ref{tab:fractional-time} refer to the case $\gamma = 0.3$ and show that \texttt{Ev-Int} significantly outperform \texttt{PGMRES}, despite the low number of iterations needed by latter to converge.
	\begin{table}[h]
		\centering
		\caption{Fractional time differentiation; performances of \texttt{PGMRES} and \texttt{Ev-Int}.}\label{tab:fractional-time}
		\begin{tabular}{|c|c|ccc|cc|ccccc|}
			\hline
			&&\multicolumn{3}{c}{\texttt{PGMRES}}&\multicolumn{2}{|c|}{\texttt{Ev-Int}}\\
			\hline
			$n$&$n_t$& Time &Res&It&Time&Res\\
			\hline
			
			\multirow{4}{*}{$1089$}& $256$ & $0.46$ & $1.0e-06$ & $5$ & $\mathbf{0.12}$ & $1.6e-09$ \\ 
			& $512$ & $0.74$ & $7.8e-05$ & $4$ & $\mathbf{0.22}$ & $1.3e-09$\\ 
			& $1024$ & $1.55$ & $2.7e-05$ & $4$ & $\mathbf{0.45}$ & $1.2e-09$\\ 
			& $2048$ & $3.68$ & $8.0e-06$ & $4$ & $\mathbf{0.87}$ & $9.9e-10$\\   
			\hline
			\multirow{4}{*}{$4225$}& $256$ & $1.56$ & $1.6e-05$ & $5$ & $\mathbf{0.39}$ & $2.3e-08$\\ 
			& $512$ & $2.77$ & $1.2e-03$ & $4$ & $\mathbf{0.87}$ & $2.0e-08$\\ 
			& $1024$ & $6.45$ & $4.1e-04$ & $4$ & $\mathbf{1.66}$ & $1.9e-08$\\ 
			& $2048$ & $14.41$ & $1.2e-04$ & $4$ & $\mathbf{3.71}$ & $1.5e-08$\\ 
			\hline
		\end{tabular}
	\end{table}
	\subsection{Wave equation}
	We consider the linear wave equation 
	$$
	\begin{cases}
		u_{tt}-\Delta u=(1+2\pi^2)e^t\sin(\pi x)\sin(\pi y),&\text{on } S\times (0, T),\\
		u(x,y,t)=0, &\text{on } \partial S \times (0, T),\\
		u(\cdot,\cdot,0) = u_0(x,y),\quad u_t(\cdot,\cdot,0) = u_1(x,y), & \text{in }S,
	\end{cases}
	$$
	with $S =(0, 1)$, $T=1$, $u_0(x,y)=u_1(x,y)=\sin(\pi x)\sin(\pi y)$. As described in \cite[Section 3]{gander2020paradiag}, by discretizing with finite differences in space and with the implicit leap-frog scheme in time we get the all-at-once linear system $(B_1\otimes I+B_2\otimes A)\mathbf x=\mathbf f$ where $A= \mathrm{trid}(-1, 2, -1)\cdot (n+1)^2$ and 
	$$
	B_1=\frac{1}{\Delta t^2}\begin{bmatrix}
		1\\
		-2&1\\
		1&-2&1\\
		&\ddots&\ddots&\ddots\\
		&&1&-2&1
	\end{bmatrix}\in\mathbb R^{n_t\times n_t},\qquad
	B_2=\frac{1}{2}\begin{bmatrix}
		1\\
		0&1\\
		1&0&1\\
		&\ddots&\ddots&\ddots\\
		&&1&0&1
	\end{bmatrix}\in\mathbb R^{n_t\times n_t}.
	$$
	We note that, for this example, all methods based on low-rank updates generate quite small Krylov subspaces, and their dimensions do not grow as the parameters $n$ and $n_t$ vary. As a consequence, \texttt{EK} is the one that achieves the cheapest consumption of computational time and we only report its performances. We remark that, in order to get comparable accuracy we have considered $d=3$ roots of the unit when running \texttt{Ev-Int}. The preconditioner employed in \texttt{PGMRES} uses the parameter $\alpha=0.1$ as suggested in \cite{gander2020paradiag} and the tolerance for stopping the GMRES iteration  has been set to $10^{-7}$. The results are shown in Table~\ref{tab:wave}. Comparing $\texttt{PGMRES},\texttt{Ev-Int}$, and $\texttt{EK}$, leads to similar conclusions as in the example of the 1D heat equation.
	\begin{table}[h]
		\centering
		\caption{Wave equation; performances of the various methods.}\label{tab:wave}
		\resizebox{1\textwidth}{!}{  
			\begin{tabular}{|c|c|ccc|cc|ccccc|}
				\hline
				&&\multicolumn{3}{c}{\texttt{PGMRES}}&\multicolumn{2}{|c|}{\texttt{Ev-Int}}&\multicolumn{5}{|c|}{\texttt{EK}}\\
				\hline
				$n$&$n_t$& Time &Res&It&Time&Res&Time&Res&rk&$T_{sylv}$ (\%)&Dim\\
				\hline
				
				\multirow{4}{*}{$1089$}& $257$ & $0.49$ & $1.2e-12$ & $3$ & $\mathbf{0.20}$ & $2.7e-10$ & $0.50$ & $2.9e-10$ & $26$ & $86.4$ & $27$                \\ 
				& $513$ & $0.72$ & $3.1e-13$ & $3$ & $\mathbf{0.34}$ & $2.6e-10$ & $0.53$ & $1.3e-10$ & $26$ & $77.5$ & $27$                \\ 
				& $1025$ & $1.42$ & $3.9e-12$ & $3$ & $\mathbf{0.71}$ & $2.5e-10$ & $\mathbf{0.71}$ & $5.3e-11$ & $27$ & $67.1$ & $27$                \\ 
				& $2049$ & $3.34$ & $1.7e-12$ & $3$ & $1.50$ & $2.5e-10$ & $\mathbf{1.25}$ & $3.7e-11$ & $27$ & $61.7$ & $27$                \\ 
				\hline
				\multirow{4}{*}{$4225$}& $257$ & $1.33$ & $3.7e-12$ & $3$ & $0.64$ & $2.8e-10$ & $\mathbf{0.57}$ & $2.9e-10$ & $26$ & $55.5$ & $27$                \\ 
				& $513$ & $2.65$ & $7.3e-12$ & $3$ & $1.27$ & $2.7e-10$ & $\mathbf{0.76}$ & $1.3e-10$ & $26$ & $54.0$ & $27$                \\ 
				& $1025$ & $5.96$ & $4.6e-11$ & $3$ & $2.96$ & $2.6e-10$ & $\mathbf{1.25}$ & $5.4e-11$ & $26$ & $34.5$ & $27$                \\ 
				& $2049$ & $12.00$ & $3.2e-11$ & $3$ & $5.38$ & $2.5e-10$ & $\mathbf{2.65}$ & $3.5e-11$ & $27$ & $34.1$ & $27$                \\ 
				\hline
		\end{tabular}}
	\end{table}
	\subsection{Runge--Kutta example}
	\upd{We present} a numerical test on the all-at-once Runge--Kutta formulation introduced in Section~\ref{sec:runge-kutta}. Let us start from the finite difference semidiscretization  of the following differential problem:
	\begin{equation}\label{eq:diff}
		\begin{cases}
			u_t=u_{xx},&(x,t)\in(0,1)^2,\\
			u(x, 0) = \sin(\pi x),
		\end{cases}\quad\Rightarrow\quad 
		\begin{cases}
			\dot U(t)=AU(t),&A=\mathrm{trid}(1,-2,1)\cdot (n+1)^2, \\
			U(0) = \mathbf x_0,&(\mathbf x_0)_j=\sin(j\pi/(n+1)).
		\end{cases}
	\end{equation}
	
	\begin{figure}[h]
		\centering
		\includegraphics[width=.5\textwidth]{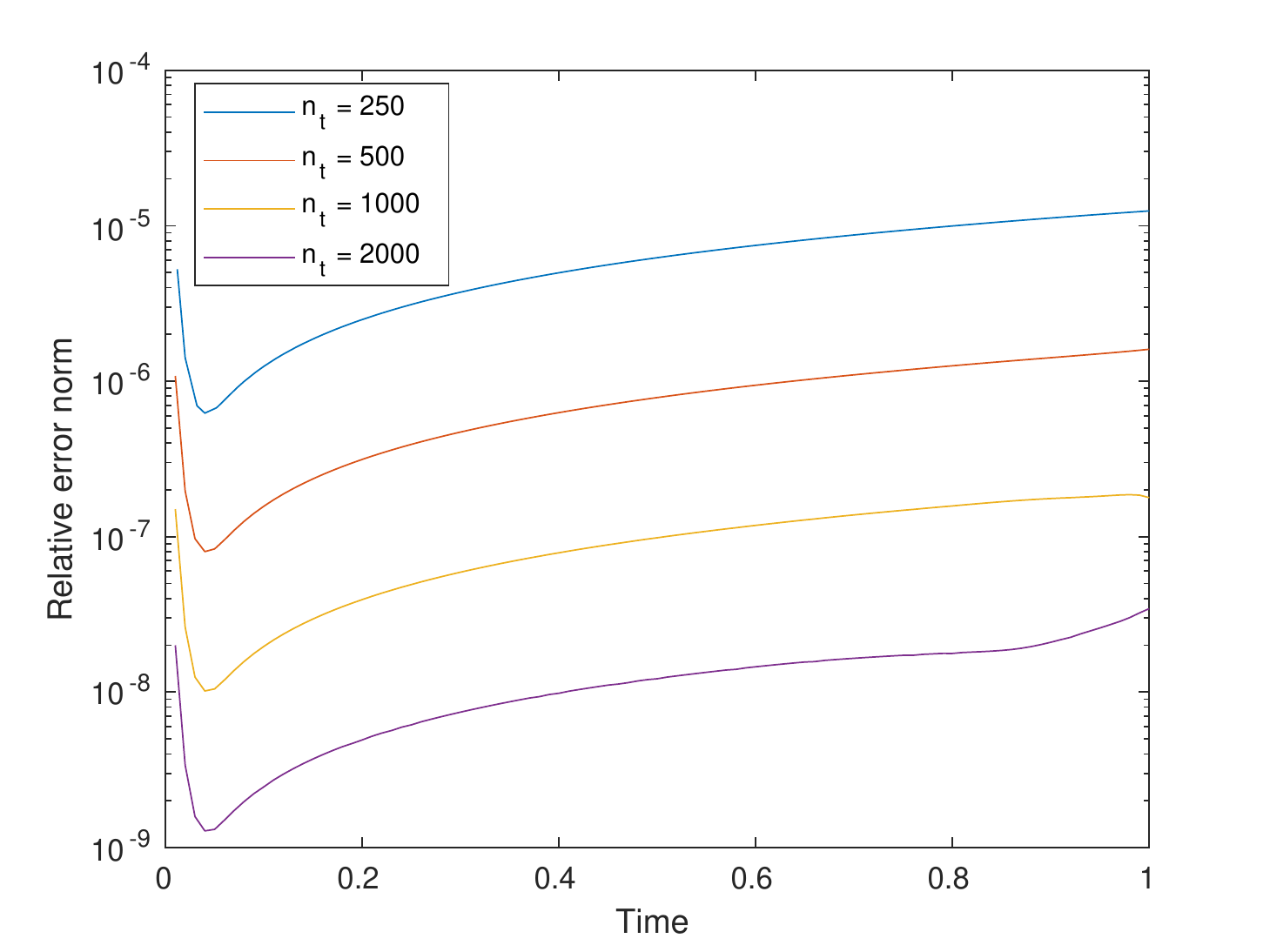}
		\caption{Relative error norm of \texttt{EK} for $n=500$ and $n_t=250,500,1000,2000$; the error decreases with order $3$ with respect to the time step size.}\label{fig:rungekutta}
	\end{figure}
	
	\begin{table}[h]
		\centering
		\caption{Performances of the Runge--Kutta schemes for different values of $n=250,500$ and $n_t=250,500,1000,2000$.}\label{tab:rungekutta}
		\resizebox{\textwidth}{!}{  
			\begin{tabular}{|c|c|ccc|cc|ccccc|}
				\hline
				&&\multicolumn{3}{c}{\texttt{PGMRES}}&\multicolumn{2}{|c|}{\texttt{Ev-Int}}&\multicolumn{5}{|c|}{\texttt{EK}}\\
				\hline
				$n(s+1)$&$n_t$& Time &Res&It&Time&Res&Time&Res&rk&$T_{sylv}$ (\%)&Dim\\
				\hline
				
				\multirow{4}{*}{$1250$}& $250$ & $0.78$ & $7.2e-08$ & $3$ & $\mathbf{0.24}$ & $1.4e-11$ & $0.36$ & $2.1e-11$ & $1$ & $69.0$ & $27$  \\ 
				& $500$ & $1.34$ & $3.2e-09$ & $3$ & $\mathbf{0.47}$ & $1.0e-11$ & $0.49$ & $1.1e-11$ & $1$ & $52.7$ & $31$  \\ 
				& $1000$ & $3.16$ & $1.5e-07$ & $3$ & $0.90$ & $1.7e-11$ & $\mathbf{0.72}$ & $2.7e-12$ & $2$ & $41.1$ & $33$  \\ 
				& $2000$ & $6.64$ & $2.8e-07$ & $3$ & $2.05$ & $5.9e-11$ & $\mathbf{1.24}$ & $4.0e-12$ & $2$ & $25.9$ & $33$  \\ 
				\hline
				\multirow{4}{*}{$2500$}& $250$ & $1.50$ & $1.8e-07$ & $3$ & $\mathbf{0.41}$ & $1.2e-10$ & $0.63$ & $3.3e-12$ & $2$ & $68.9$ & $27$  \\ 
				& $500$ & $3.21$ & $1.6e-07$ & $3$ & $\mathbf{0.83}$ & $6.6e-11$ & $0.88$ & $4.9e-12$ & $2$ & $53.5$ & $31$  \\ 
				& $1000$ & $6.49$ & $3.4e-06$ & $3$ & $1.80$ & $6.2e-11$ & $\mathbf{1.39}$ & $7.2e-12$ & $2$ & $39.7$ & $33$  \\ 
				& $2000$ & $13.35$ & $7.9e-05$ & $3$ & $3.73$ & $2.2e-10$ & $\mathbf{2.41}$ & $1.1e-11$ & $2$ & $24.3$ & $33$  \\ 
				\hline
		\end{tabular}}
	\end{table}
	\noindent We consider solving the all-at-once formulation \eqref{eq:runge-kutta-lin-sys} corresponding to the  following  four-stage ($s=4$), $3$rd order, \upd{\emph{Diagonally Implicit Runge–Kutta method} (DIRK)}:
	$$
	\begin{array}{c|cccc}
		1/2&1/2\\
		2/3&1/6&1/2\\
		1/2&-1/2&1/2&1/2\\
		1&3/2&-3/2&1/2&1/2\\
		\hline
		&3/2&-3/2&1/2&1/2
	\end{array}.
	$$
	More specifically, we solve \eqref{eq:runge-kutta-lin-sys} by applying the methods \texttt{PGMRES}, \texttt{Ev-Int} and \texttt{EK} adapted as described in Section~\ref{sec:runge-kutta}. For this experiment, the parameter $\rho$ used by \texttt{Ev-Int} is set to the value $0.05$.  As shift parameter for the update equation \eqref{eq:lowrank-sylv} solved by \texttt{EK}, it has been chosen $\sigma = -2$.
	The performances of the $3$ methods are reported in Table~\ref{tab:rungekutta}; as observed in previous tests, the approaches based on low-rank update and evaluation interpolation outperforms \texttt{PGMRES} for both speed and accuracy. As the number of time steps increases, the cost of the update equation becomes relatively cheap, making \texttt{EK} faster than \texttt{Ev-Int}. The relative residuals obtained with \texttt{EK} and \texttt{Ev-Int} are comparable.
	
	Problem~\eqref{eq:diff} admits the explicit solution $U(t)=\exp(tA)\mathbf x_0$ so that we can compare approximate solutions with respect to a benchmark solution computed by means of the \texttt{expm} Matlab function. Figure~\ref{fig:rungekutta} reports the relative error $\norm{\widehat {\mathbf x}_t-\exp(tA)\mathbf x_0}_2/\norm{\exp(tA)\mathbf x_0}_2$ of the approximants $\widehat{\mathbf x_t}$ returned by \texttt{EK}, for $n=500$ and different time step sizes; the plot is in agreement with the order $3$ of the Runge--Kutta scheme.
	\upd{
		\subsection{Preliminary experiments on parallelism}\label{sec:parallel-exp}
		We conclude by testing the parallel scaling features of the proposed algorithms in a multi-core computational environment. The experiment in this section has been run on a server with two Intel(R) Xeon(R) E5-2650v4 CPU with
		12 cores and 24 threads each, running at 2.20 GHz, using MATLAB R2022b with the Intel(R)
		Math Kernel Library Version 11.3.1. The test has been run using the SLURM scheduler,
		allocating 12 cores and 250 GB of RAM. Our implementation exploits parallelism for the solution of the block diagonal linear system generated by \textsc{fast\_diag\_solve} (lines 4-6 of Algorithm~\ref{alg:fast-diag}). This parallel routine is used as a preconditioner in \texttt{PGMRES}, it is called at line~\ref{step:X0} of Algorithm~\ref{alg:update}, and it is called at line~4 of \texttt{Ev-Int}. Finally, we also parallelize the for loop in \texttt{Ev-Int} (lines 2-5 of Algorithm~\ref{alg:ev-int}).
		
		As case study, we consider the fractional diffusion example of Section~\ref{sec:frac-space} with $n=4225$ and $n_t=8192$, and  we measure the execution time of the procedures \texttt{PGMRES}, \texttt{Ev-Int}, and \texttt{ZOL-DI}, as the number of cores $p$ ranges in the set $\{1,2,4,8,12\}$. The tolerance to stop GMRES has been set to $10^{-6}$ in order to get comparable residual norms with the other two methods.
		To evaluate the efficiency of the implementation, we also compute the strong scaling efficiency (SE) with $p$ cores defined as
		$$
		\mathrm{SE}(p):=100 \cdot \frac{\mathrm{Time}(\text{$1$ core})}{p\cdot\mathrm{Time}(\text{$p$ cores})}.
		$$
		The results reported in Table~\ref{tab:parallel} clearly show that \texttt{Ev-Int} is the method that gains the most from parallelism, maintaining more than $70$\% of strong scaling efficiency, up to $12$ cores. Also \texttt{PGMRES} leverages the use of  multiple cores to the extent of becoming faster than \texttt{ZOL-DI}, when $p=12$. On the contrary, our current implementation of \texttt{ZOL-DI} does not exploit more than one core for solving the update equation and, consequently, parallel efficiency drops as soon as the cost of \eqref{eq:update-eq} becomes dominant. An adaptation of the parallel implementation of RKSM~\cite{Berljafa2017} is nontrival and beyond the scope of this work, but it potentially addresses the scaling issue of \texttt{ZOL-DI}. 
		
		\begin{table}[h]
			\centering
			\caption{Fractional space diffusion with $n=4225$ and $n_t=8192$; parallel efficiency of the procedures as the number of cores $p$ increases.}\label{tab:parallel}
				\begin{tabular}{|c|cc|cc|ccc|}
					\hline
					&\multicolumn{2}{c}{\texttt{PGMRES}}&\multicolumn{2}{|c|}{\texttt{Ev-Int}}&\multicolumn{3}{|c|}{\texttt{ZOL-DI}}\\
					\hline
					$p$& Time &SE&Time&SE&Time&$T_{sylv}$ (\%)&SE \\
					\hline
					$1$&$393.23$&$100.0$&$128.18$&$100.0$&$\mathbf{118.27}$&$46.3$&$100.0$
					\\ 
					$2$&$206.50$&$95.2$&$\mathbf{66.91}$&$95.8$&$87.30$&$62.5$&$67.7$
					\\ 
					$4$&$135.94$&$72.3$&$\mathbf{38.13}$&$84.0$&$71.45$&$74.5$&$41.4$
					\\ 
					$8$&$74.44$&$66.0$&$\mathbf{20.66}$&$77.5$&$63.02$&$84.9$&$23.5$
					\\ 
					$12$&$53.64$&$61.1$&$\mathbf{14.91}$&$71.6$&$60.18$&$89.2$&$16.4$
					\\ 
					\hline
				\end{tabular}
		\end{table}

	}

	\section{Conclusions}
	We  have  proposed  a  tensorized Krylov  subspace  method and an interpolation method  for replacing the iterative refinement phase that is  required by ParaDiag algorithms with uniform time steps. In the case of multistep methods,  a theoretical analysis of the approximation property of the Krylov subspace method  has  been  provided. We have shown several  numerical  tests demonstrating that our  approaches  can  significantly  outperform  block-circulant preconditioned GMRES iterations for solving linear initial value problems. Let us emphasize that the focus of this work was on theoretical and algorithmic aspects; we have not considered the concrete parallel implementation of our algorithms, which will be subject to future work. From the analysis in Sections~\ref{sec:cost} and~\ref{sec:costev} we expect that the interpolation method (Algorithm~\ref{alg:ev-int}) will perform significantly better in a  massively parallel environment.
	
	We remark that Algorithm~\ref{alg:update} and Algorithm~\ref{alg:ev-int} can directly be put in place instead of the existing ParaDiag algorithms for solving certain non-linear ODEs of the form $M\dot U+f(U)=0$ \cite{gander2017} and for computing the so-called \emph{Coarse Grid Correction} of the Parareal algorithm~\cite{lions2001,gander2007}.
	
	\upd{In Section~\ref{sec:runge-kutta},} we have discussed how to extend the methodology to the case of \upd{implicit} Runge--Kutta time integration, but some questions remain for further work. In particular, a careful analysis of the low-rank approximability of equation~\eqref{eq:lowrank-sylv} and of the choice of the shift parameter $\sigma$, would be desirable.
	
	\begin{paragraph}{Acknowledgments.}
		The authors thank Bernd Beckermann for inspiring discussions on the Zolotarev problems in Section~\ref{sec:zolotarev}. 
	\end{paragraph}
	\bibliographystyle{abbrv}
	\bibliography{Bib} 
	
\end{document}